\documentclass[a4paper,11pt,twoside]{amsart}
\input{xy}
\xyoption{all}
\usepackage{master}
\begin{document}
\bibliographystyle{hamsalpha}
\comment{
\title[...]{...}                                 
    \author[...]{...}                                
    \address{...}                                    
    \curraddr{...}                                   
    \email{...}                                      
    \urladdr{...}                                    
    \dedicatory{...}                                 
    \date{...}                                       
    \thanks{...}                                     
    \translator{...}                                 
    \subjclass{...}           
}
\title[Vector bundles on Mumford curves]{On representations attached to semistable vector bundles on Mumford curves}
\author{Gabriel Herz}
\address{Gabriel Herz\\ Mathematisches Institut der Universit\"at M\"unster\\ Einsteinstra\ss{}e 62\\ 48149 M\"unster\\ Germany}
\email{gherz@math.uni-muenster.de}
\thanks{This is a slightly modified version of my PhD-thesis at the University of M\"unster under supervision of Christopher Deninger.}
\date{\today}
\maketitle
\tableofcontents
\section*{Introduction}
A classical result by Andr\'e Weil, proved in 1938, asserts that a holomorphic vector bundle on a Riemann surface is given by a representation of the fundamental group if and only if each indecomposable component is of degree zero (cf. \cite{Weil1938}). Furthermore it is known that unitary representations of the fundamental group are in one-to-one correspondance to polystable vector bundles of degree 0 (cf. \cite{Nara1965}).\\

In 1983 Gerd Faltings introduced the notion of {\em $\phi$-bounded} representations and proved a corresponding result in $p$-adic analysis for vector bundles on a Mumford curve over a discrete non-Archimedean field. He proved an equi\-valence of categories between semistable vector bundles of degree zero on a Mumford curve and $\phi$-bounded representations of its Schottky group. In his proof he used the theory of formal schemes and was therefore limited to discrete fields. In 1986 Marius van der Put and Marc Reversat generalised Faltings' result to arbitrary non-Archimedean fields by using methods from rigid geometry. Unfortunately the functor they constructed does not commute with duals or tensor products. In the case of line bundles it is easy to see that this is inherently caused by the definiton of $\phi$-boundedness. In the sequel we will call their representations {\em PR-representations}.\\
Let $X$ be a projective, smooth and geometrically connected algebraic curve over a finite field extension of $\Qp$. In 2004 Annette Werner and Christopher Deninger constructed an \'etale parallel transport for vector bundles with {\em potentially strongly semistable reduction} on $X_{\Cp}$ (cf. \cite{Deni2004}). Restricted to the algebraic fundamental group this is a functor that associates a continuous $\Cp$-vector space representation (in the following called {\em DW-representation}) of the algebraic fundamental group of $X_{\Cp}$ to every vector bundle of this class. This functor is $\Cp$-linear, exact and commutes with duals, tensor products, internal homs and exterior powers. On can pose the question, whether the category of vector bundles with potentially strongly semistable reduction is equal to the category of semistable vector bundles of degree zero. In this thesis we will only be concerned with a special case of the Deninger-Werner construction.\\
Recently Gerd Faltings had announced a $p$-adic version of non-abelian Hodge-theory. In \cite{Falt2003} he defines a category of certain \emph{generalised representations} and proves an equivalence of categories between them and vector bundles on $X_{\Cp}$ endowed with a $p$-adic Higgs field. In his proof he uses the theory of almost \'etale extensions. The representations of $\pi^{alg}_{1}(X_{\overline{\Qp}},x)$ form a full subcategory of Faltings' generalised representations and Faltings suggests that under the equivalence of categories semistable vector bundles of degree zero come from representations of $\pi^{alg}_{1}(X_{\overline{\Qp}},x)$. In this thesis we will not rely on the greater generality of Faltings approach.

Let $X$ be a Mumford curve over a finite field extension of $\Qp$. In this thesis we compare the DW-representations attached to a class of semistable vector bundles of degree zero on $X_{\Cp}$ to the PR-representations defined for this class of vector bundles. Performing this, obvious problems occur:
\begin{enumerate}
\item it is not known whether the DW-representation does exist for all semistable vector bundles of degree zero.
\item Different groups are represented.
\item The functors of Deninger--Werner and van der Put--Reversat have different properties. For example the DW-functor commutes with duals and tensor products, the PR-functor in general does not.
\end{enumerate}
The solution of all of these problems is 
\begin{enumerate}
\item to consider only the subcategory of semistable vector bundles of degree 0 that have a vector bundle model on the extension of the minimal regular model of $X$ to the ring of integers $\op$ of $\Cp$, 
\item to introduce the notion of finite topological coverings and the finite topological fundamental group of $X$ 
\item and at last to prove that the DW-representation factorises through this fundamental group.
\end{enumerate} 
Having done this, we can prove that the DW-representation and the pro-finitely completed PR-representation attached to vector bundles in the previously mentioned class are isomorphic. The PR-representations attached to vector bundles in this class are isomorphic to representations that have image in $\Gl_{\rk}(\op)$ if $\rk$ is the rank of the vector bundle considered.\\
At least for line bundles the previously mentioned class is the best possible on which both representations agree, since for line bundles whose PR-representa{\-}ti{\-}on is not represented by numbers of norm equal to one, the PR-represen{\-}ta{\-}tion does not commute with tensor products and duals, but the DW-functor does.\\
Introducing the fundamental group of finite topological coverings suggests itself, since we have to compare representations of the topological fundamental group (which is the Schottky group) with representations of the finite \'etale fundamental group. The coverings that are topological and finite \'etale are exactly the finite topological coverings.\\ 
We prove our result by extensive use of GAGA theorems between rigid, formal and algebraic geometry. Because we have some GAGA results only in the case of discrete valuation we reduce the case that the vector bundle is only defined after extension to $\Cp$ to the case that it is already defined over a discrete valuation field by an argument from non-abelian cohomology.\\ 
As an illustration we will have a closer look at Mumford curves of genus 1 and 2 and at vector bundles on them. \\

This thesis is organised as follows. In the first section we remind the reader of some notions of rigid and analytic spaces, which are important for us. We cite some GAGA results and prove slight extensions of them. Especially quotients of schemes by finite groups are considered more closely. We state the basic notions of Galois theory, introduce the finite topological fundamental group and prove that it satisfies the six axioms of Grothendieck for a Galois theory. At the end of the first section we introduce Mumford curves and discuss their stable and minimal regular model.\\
In the second section we describe the constructions of van der Put and Reversat, Faltings and  Deninger--Werner. We characterise the semistable vector bundles of degree zero that have a vector bundle model on the minimal reguar model of the Mumford curve by their PR-representation and we compare the construction of Faltings with the one of van der Put--Reversat. In the next subsection we compare the DW-representation with the PR-representation. We deduce the case that the vector bundle is only defined after base change to $\Cp$ from the case that it is already defined over a discrete valuation field and prove this case first.\\
In the last section we give some illustrations. We study the various Galois groups of a Tate curve and investigate vector bundles on Tate curves and on Mumford curves of genus 2 in more detail.\\
I would like to thank my supervisior Prof. Christopher Deninger and Prof. Annette Werner for introducing my to this interesting topic. The final version benefited from discussions with Jan Kohlhaase, Sylvain Maugeais, Roland Olbricht, Matthias Strauch and Stefan Wiech. I am grateful that Sylvain Maugeais and Stefan Wiech read a preliminary version.\\
This thesis was partially supported by the {\em Deutsche Forschungsgemeinschaft} at the SFB 478 {\em Geometrische Strukturen in der Mathematik}.
\section{Preliminaries}\label{chpre}
\subsection{Notations and conventions}
We use the theories of schemes, formal schemes, rigid spaces and Berkovich spaces and assume that the reader is familiar with these theories. In particular for the theory of rigid spaces the reader might consult \cite{Bosc1984}, for Berkovich spaces \cite{Berk1990} and \cite{Berk1993}, for formal schemes the reader might refer to \cite{EGAn}.\\
A field endowed with a non-Archimedean valuation that is complete with respect to this valuation is called {\em non-Archimedean field}. A {\em discrete non-Archimedean field} is a non-Archimedean field whose valuation is discrete. We will assume all non-Archimedean fields to have non-trivial valuation and to be of characteristic zero. \\ 
The following notations will often be used.\\
The term $\M(A)$ denotes the {\em Berkovich spectrum} of a commutative Banach ring $A$ with unit, that is it denotes the set of all bounded multiplicative semi-norms on $A$ provided with the weakest topology with respect to which all real-valued functions of the form $|\cdot|\mapsto |f|$ $(f\in A)$ are continuous. Let $\K$ be a non-Archimedean field. The terms {\em $\K$-analytic space} and {\em strictly $\K$-analytic space} denote the Berkovich analytic spaces that are defined in \cite[page 22]{Berk1993}. A $\K$-analytic space is called {\em good} if every point has an affinoid neighbourhood. The term {\em analytic space} always means {\em Berkovich analytic space}.\\
A Hausdorff topological space is called {\em paracompact} if every open covering has a locally finite refinment. A rigid space is called {\em quasiseperated} if the intersection of two open affinoid domains is the finite union of open affinoid domains.\\
If $f$ is an element of a $\K$-affinoid algebra $A$, then we denote by $\spn{f}$ the spectral norm of $f$, which is defined as the supremum of all $\n{f(x)}$, where $x$ runs through the set of all maximal ideals of $A$.\\ 
We denote the formal spectrum of an adic ring by $\Spf(\cdot)$, the maximal spectrum of an affinoid algebra by $\Sp(\cdot)$, the spectrum of maximal ideals of a scheme by $\Spm(\cdot)$ and the ordinary spectrum of a scheme by $\Spec(\cdot)$.\\
If $k$ is a field, then the term {\em algebraic $k$-variety} denotes a geometrically integral and seperated scheme of finite type over $k$. The term {\em $k$-curve} denotes an algebraic $k$-variety of dimension 1.\\
Let $\K$ be a non-Archimedean field, $\K^{\circ}$ its ring of integers and $\K^{\circ\circ}$ its unique maximal ideal. If in particular $\K=\Cp$, then we also denote $\K^{\circ}$ by $\op$. Let $\L$ be a non-Archimedean field, too. In the following, when we write {\em $\K$ is a subfield of $\L$} resp. {\em $\K \subset \L$} we assume that the norm of $\K$ is induced by the norm of $\L$.  
We call a finite extension field of $\Qp$ a {\em local number field}. Local number fields are complete with respect to the norm that is induced by $\Qp$.\\
Let $R$ be a valuation ring with quotient field $Q$ and $X$ a $Q$-curve. A finitely presented, proper and flat $R$-scheme $\X$ with the property $\X\otimes_{R} Q=X$ is called {\em $R$-model} of $X$.
If $R$ is a Dedekind domain, then we call an integral, projective and flat R-scheme of dimension 2 a {\em projective flat R-curve}.\\
Following Mumford \cite{Mumf1962}, a vector bundle $E$ on an algebraic $\K$-curve is {\em stable} if for all proper subbundles $F$,
\[\frac{\deg F}{\deg E}<\frac{\rank F}{\rank E}.\] The vector bundle $E$ is called {\em semistable} if for all subbundles $F$,
\[\frac{\deg F}{\deg E}\leq\frac{\rank F}{\rank E}.\] 
\subsection{Reduction and formal rigid spaces}
In this subsection we explain reduction of rigid spaces and define the notion of a formal rigid space. Then we compare formal rigid spaces to formal schemes and models of algebraic curves to their associated formal rigid spaces. 
\subsubsection{Reduction of rigid spaces}\text{ }\\
Let $\K$ be a non-Archimedean field. To define the reduction of rigid spaces we have to define the reduction of affinoids first.\\
Let $A$ be an affinoid $\K$-algebra and denote by $||\cdot||_{sp}$ the spectral norm; then \[A^{\circ}:=\{a\in A\suchthat\spn{a}\leq 1\}\] is a ring, and \[A^{\circ\circ}:=\{a\in A \suchthat \spn{a} < 1\}\] is an ideal in $A$. The residue ring $A^{\circ}/A^{\circ\circ}$ is denoted by $\widetilde{A}$.
Let $\Sp(A)$ denote the maximal spectrum of $A$. For every point $x\in\Sp(A)$ there is a map $\phi_{x}:A\ra A/\m_{x}$, where $\m_{x}$ is the maximal ideal associated to $x$. If we denote $A/\m_{x}$ by $\L_{x}$, then the field $\L_{x}$ is a finite extension of $\K$, hence it carries a unique extension of the norm on $\K$. Therefore there is an induced map $\widetilde{\phi_{x}}:\widetilde{A}\ra\widetilde{\L_{x}}$. The set $\Kern\widetilde{\phi_{x}}$ is a maximal ideal of $\widetilde{A}$. So we get a map 
\[\pi:\Sp(A)\ra\Spm(\widetilde{A}), x\mapsto\Kern(\widetilde{\phi_{x}}).\] 
\begin{defin}
\begin{enumerate}
\item For $X:=\Sp A$ and $\widetilde{X}^{c}:=\Spm(\widetilde{A})$ this map is called {\em the canonical reduction map} 
\[\Red^{c}_{X}:X\ra\widetilde{X}^{c}.\]
\item An affinoid domain $V$ in $X$ is said to be {\em formal} if the induced morphism of the canonical reductions $\widetilde{V}^{c}\ra\widetilde{X}^{c}$ is an open immersion.
\end{enumerate}
\end{defin}
\begin{nota}
Let $X$ be a $\K$-scheme of finite type and denote its topological subspace of closed points by $i:X_{\circ}\inj X$. The scheme $(X_{\circ},i^{-1}\Loc_{X})$ is dense in $(X,\Loc_{X})$. Therefore it is no drawback to have the reduction map only to the maximal spectrum of $A$, since $\Spm A=(\Spec A)_{\circ}$.  
\end{nota}
It is always the case that the reduction preimage of a Zariski open subset of $X$ is admissible in the weak Grothendieck topology. But one can strengthen this by some conditions as in the following lemma.  
\begin{lemma}[\cite{Fres2004}, Lemma 4.8.1]\label{coverred}
Suppose $X:=\Sp A$ is reduced and $\spn{A}\subset\n{\K}$. If $U$ is an open affine subset of $\widetilde{X}^{c}$, then
\begin{enumerate}
\item $(\Red^{c}_{X})^{-1}(U)$ is affinoid,
\item the ring $\Loc((\Red^{c}_{X})^{-1}(U))$ is reduced and its spectral norm takes values in $\n{\K}$, and
\item its reduction $\widetilde{\Loc((\Red^{c}_{X})^{-1}(U))}$ is canonically isomorphic to the affine $\widetilde{\K}$-algebra of the regular functions on $U$.
\end{enumerate}
\end{lemma}

\begin{defin}[\cite{Fres2004}, Definitions 4.8.3 ]\label{reddef}
Let $X$ be a reduced rigid space.
\begin{enumerate}
\item An admissible affinoid covering $(U_{i})_{i\in I}$ of $X$ is called {\em formal covering} or {\em pure covering} if the following holds:
\begin{enumerate}
\item the set $\spn{\Loc(U_{i})}$ is a subset of $\n{\K}$,
\item for all $i\in I$ there are only finitely many $j\in I$ such that $U_{i}\cap U_{j}\ne \emptyset$,
\item \label{redaffine} if $U_{i}\cap U_{j}\ne \emptyset$, then $U_{i}\cap U_{j}$ is formal open, and
\item the natural map $\Loc(U_{i})^{\circ}\widehat{\otimes}_{\K^{\circ}}\Loc(U_{j})^{\circ}\ra\Loc(U_{i}\cap U_{j})^{\circ}$ is surjective.
\end{enumerate}
\item By glueing affine varieties one obtains a global reduction map 
\[\Red:=\Red_{X,(U_{i})_{i\in I}}:X\ra \widetilde{(X,(U_{i})_{i\in I})}.\]
\item Two coverings $(U_{i})_{i\in I}$ and $(V_{j})_{j\in J}$ are said to be {\em equivalent} if $U_{i}\cap V_{j}$ is a formal open subset of both $U_{i}$ and $V_{j}$ for all $i\in I$ and $j\in J$.   
\item Let $\U$ be an equivalence class of formal coverings. A pair $(X,\U)$ is called {\em formal rigid $\K$-space}. We denote by $(X,(U_{i})_{i\in I})$ a formal rigid $\K$-space whose equivalence class of formal coverings is represented by $(U_{i})_{i\in I}$. The reduction of a formal rigid space does only depend on the eqiuvalence class of its formal covering.   \item A {\em morphism of formal rigid spaces} $\phi:(X,\U)\ra (Y,\V)$ is a morphism of rigid spaces $\phi:X\ra Y$ with the additional property that for every formal open subset $V$ the preimage $\phi^{-1}(V)$ is formal open, too. That means that 
\[\Red_{X,(U_{i})_{i\in I}}(\phi^{-1}(V))\ra\Red_{X,(U_{i})_{i\in I}}(X)\] is an open immersion.
\end{enumerate}
\end{defin}
We want to prove an equivalence of categories between formal rigid spaces and a certain class of admissible formal schemes.
\begin{defin}\text{ }
\begin{enumerate}
\item Let $R$ be an $I$-adic ring. A topological $R$-algebra $A$ is called {\em admissible} if it has no $I$-torsion and if there is a finitely generated ideal $J$ in the ring of restricted power series $R\{T_{1},\ldots,T_{n}\}$ such that $A$ is isomorphic to an $R$-algebra of type $R\{T_{1},\ldots,T_{n}\}/J$ endowed with the $I$-adic topology. 
\item A formal $R$-scheme is {\em admissible} if and only if there is an open affine covering $(\Spf A_{i})_{i\in I}$ such that each $A_{i}$ is an admissible $R$-algebra.
\end{enumerate}
\end{defin}
The following lemma shows that every morphism between affinoid $\K$-algebras is contractive. Its proof is a slight variation of a proof given in a lecture by S. Bosch.
\begin{lemma}\label{contractive}
If $\phi:A\ra B$ is a morphism between affinoid $\K$-algebras, then $\spn{\phi(a)}\leq\spn{a}$ for all $a\in A$. By definiton, this means that $\phi$ is {\em contractive}.
\end{lemma}
\begin{proof}
Assume the contrary, that is there is a function $a$ in $A$ and a maximal ideal $\m_{B}$ in $B$ such that for all maximal ideals $\m_{A}$ in $A$ it is true that $\n{a(\m_{A})}<\n{\phi(a)(\m_{B})}$. We prove that this assertion is false. Define $\m:=\phi^{-1}(\m_{B})$, this is a maximal ideal. Thus for all $a\in A$ it is true that $\n{a(\m)}=\n{a(\phi^{-1}(\m_{B}))}=\n{\phi(a)(\m_{B})}$.
\end{proof}
Now we can state and prove the announced equivalence.
\begin{lemma}
Let $\K$ be a non-Archimedean field and $\K^{\circ}$ its ring of integers. There is an equivalence between the category of formal rigid $\K$-spaces $(X,(U_{i})_{i\in I})$ and the category of separated and admissible formal $\K^{\circ}$-schemes $\X$ that admit a locally finite covering. The underlying topological spaces of $\widetilde{(X,(U_{i})_{i\in I})}$ and $\X$ coincide.
\end{lemma}
\begin{proof}
The existence of the functor $\cdot^{rig}$ from formal spaces to rigid spaces is well known. On affine spaces the functor maps \[\Spf(A)\mapsto \Spm(A\otimes \K)\] for an admissible $\K^{\circ}$-algebra $A$. It is known that $\cdot^{rig}$ is essentially surjective. For these facts see \cite{Rayn1974}.\\
Let $(X,(U_{i})_{i\in I})$ be a formal rigid space, in particular there are affinoid $\K$-algebras $A_{i}$ such that $U_{i}:=\Sp(A_{i})$. Define $\X_{i}:=\Spf(A_{i}^{\circ})$, the family $(\X_{i})_{i\in I}$ is a locally finite cover of the formal scheme $\X$. The formal scheme $\X$ does only depend on the equivalence class of $(U_{i})_{i\in I}$, it is separated by properties c) and d) of a formal rigid cover and it is admissible.\\
Let $\X$ be a separated and admissible formal $\K^{\circ}$-scheme, and let $(\Spf(A_{i}))_{i\in I}$ be a locally finite cover of $\X$. The family $(\Sp(A_{i}\otimes\K))_{i\in I}$ is a formal cover of $X:=\X^{rig}$.
The coincidence of $\widetilde{(X,(U_{i})_{i\in I})}$ and $\X$ as topological spaces follows by lemma \ref{coverred}.\\
Let $f:(X,(U_{i})_{i\in I})\ra (Y,(V_{j})_{j\in J})$ be a morphism between formal rigid $\K$-spaces. We have to construct a morphism between the formal $\K^{\circ}$-schemes $\widetilde{f}:\X\ra\Y$. It is enough to construct $\widetilde{f}$ locally. By definition of a morphism between formal rigid spaces, the preimage of a formal open set is formal open. Because $Y$ is covered by formal open sets, and formal open sets are affinoid (lemma \ref{coverred}), it is enough to construct $\widetilde{f}$ if $f$ is a morphism between affinoids.\\
Given affinoids $U$ and $V$, the existence of the morphism $f:U\ra V$ is equivalent to the existence of $\Loc(V)\ra\Loc(U)$, equivalent to the existence of the morphism $\Loc(V)^{\circ} \ra\Loc(U)^{\circ}$, because $\Loc(V)\ra\Loc(U)$ is contractive by lemma \ref{contractive}, equivalent to the existence of $\Spf(\Loc(U)^{\circ})\ra\Spf(\Loc(V)^{\circ})$.\\ 
If on the other hand $f:\X\ra\Y$ is a morphism between separated and admissible formal schemes, then we can choose locally finite coverings $(U'_{i})_{i\in I}$ of $\X$ and $(V'_{j})_{j\in J}$ of $\Y$ such that for every $i\in I$ there is a $j\in J$ with the property $f^{-1}(V'_{j})=U'_{i}$. This is equivalent to the existence of a morphism $\Loc_{V'_{j}}\ra\Loc_{U'_{i}}$. The existence of this morphism is equivalent to the existence of $\Loc_{V'_{j}}\otimes\K\ra\Loc_{U'_{i}}\otimes\K$, which is equivalent to the existence of a morphism $\Sp(\Loc_{U'_{i}}\otimes\K)\ra\Sp(\Loc_{V'_{j}}\otimes\K)$. This defines a morphism between the rigid spaces that are associated to $\X$ and $\Y$.    
\end{proof}
Projective models of algebraic curves can be classified by formal rigid structures.
\begin{theorem}[\cite{Fres2004}, Theorem 4.10.7]\label{curvesdvr}
Let $\K$ be a discrete non-Archi\-me\-dean field and $X$ a regular and projective curve over $\K$. Then there is a bijection between
\begin{enumerate}
\item formal rigid structures $(X^{rig},(U_{i})_{i \in I})$ on the rigidification of $X$, and
\item projective $\Ok$-models of $X$ whose reduced special fibre is equal to $\widetilde{(X^{rig},(U_{i})_{i \in I})}$.
\end{enumerate}
For every analytic reduction $\widetilde{(X^{rig},(U_{i})_{\in I})}$ of $X^{rig}$ the genus of $X$ is equal to the arithmetic genus of $\widetilde{(X^{rig},(U_{i})_{\in I})}$. Moreover, the formal scheme associated to $(U_{i})_{i \in I}$ coincides with the formal completion of the corresponding $\Ok$-scheme $\X$ at its closed fibre.
\end{theorem}

\subsubsection{Formal analytic vector bundles}\text{ }\\
Instead of writing $(X,(U_{i})_{i\in I})$ for a formal $\K$-rigid space one may also write $(X,r,\widetilde{X})$, where $\widetilde{X}$ is the specific reduction associated to $(U_{i})_{i\in I}$ and $r$ denotes the reduction map $r:X\ra\widetilde{X}$.
\begin{defin}[\cite{Fres2004}, page 186]\text{ }
\begin{enumerate}
\item A {\em formal vector bundle} $E$ of rank $\rk$ on $(X,r,\widetilde{X})$ is a sheaf of $r_{*}\Loc_{X}^{\circ}$-modules on $\widetilde{X}$ that is locally isomorphic to $(r_{*}\Loc_{X}^{\circ})^{\rk}$.
\item On $\widetilde{X}$ the vector bundle $E\otimes\widetilde{\K}$ is defined by \[(E\otimes\widetilde{\K})(U):=E(U)\otimes_{\K^{\circ}}\widetilde{\K}.\]
\item We write $E\otimes\K$ for the unique rigid vector bundle on $X$ that satisfies \[(E\otimes\K)(r^{-1}(U))=E(U)\otimes_{\K^{\circ}}\K\] for every open affine $U\subset\widetilde{X}$.
\end{enumerate}
\end{defin}
\begin{nota}
Denote by $\widehat{\X}$ the formal completion of $\X$. It is true that $r_{*}\Loc_{X}^{\circ}=\Loc_{\fX}$, and the underlying topological spaces of $\fX$ and $\widetilde{X}$ are equal; therefore every formal rigid vector bundle is an $\Loc_{\fX}$-vector bundle, hence a formal vector bundle and vice versa. By the identification \[r_{*}\Loc_{X}^{\circ}=\Loc_{\fX}\] 
the category of formal rigid vector bundles is equivalent to the category of formal (algebraic) vector bundles.
\end{nota}
\subsection{GAGA-results}
In this subsection we collect various results of GAGA type. We compare sheaves and vector bundles on different sites and investigate the GAGA of quotients by finite groups. For the sake of completeness we cite the relevant theorems in the first paragraph. 
\subsubsection{Rigidification, analytification and algebraisation}\text{ }\\
Let $\K$ be a non-Archimedean field.
Rigidification of algebraic varieties is presented in \cite{Fres2004}:
let $\text{Alg}_{\K}$ be the category of algebraic $\K$-varieties. Let $\text{Rigid}_{\K}$ be the category of rigid $\K$-spaces provided with the weak Gro\-then\-dieck topology. There is a functor \[\text{{\rm rig}}:\text{{\rm Alg}}_{\K}\ra\text{{\rm Rigid}}_{\K}\text{, }X\mapsto X^{rig}.\] It commutes with fibre products and $X^{rig}$ is separated.
For this statement see for example \cite[Example 4.3.3]{Fres2004}.\\
The analogous result for Berkovich spaces is the following.
Let $X$ be a $\K$-scheme of locally finite type and $\Hom(\cdot,\cdot)$ be the set of homomorphisms in the category of $\K$-ringed spaces. Define the functor 
\[\Phi:Y\mapsto \Hom_{K}(Y,X)\] from the category of good $\K$-analytic spaces to the category of sets.
\begin{theorem}[\cite{Berk1990}, Theorem 3.4.1, Corollaries 3.4.14, 3.4.12, 3.4.13]\label{curveana}
The functor $\Phi$ is represented by a closed $\K$-analytic space $X^{an}$ with a morphism $\pi:X^{an}\ra X$. They have the following properties:
\begin{enumerate}
\item There is a bijection $X^{an}(\L)\ra X(\L)$ for any non-Archimedean field $\L/\K$. Furthermore, the map $\pi$ is surjective and induces a bijection $X_{\circ}^{an}\ra X_{\circ}$. Here $X^{an}_{\circ}$ resp. $X_{\circ}$ denotes the set of closed points $\{x\in X | [\HH(x):\K]<\infty\}$.
\item On proper $\K$-schemes the functor $\cdot^{an}$ is a fully faithful functor.
\item Given a reduced proper $\K$-scheme $X$, the functor $\cdot^{an}$ gives an equivalence of categories between the finite (resp. finite \'etale) schemes over $X$ and the finite (resp. finite \'etale) good $\K$-analytic spaces over $X^{an}$. 
\item Every reduced proper good $\K$-analytic space $X$ of dimension 1 is algebraic; i.e. there is a projective algebraic $\K$-curve $Y$ such that $X\cong Y^{an}$.
\end{enumerate}
\end{theorem}
The following theorem by Berkovich gives a GAGA result for rigid and analytic spaces.
\begin{theorem}[\cite{Berk1993},Theorem 1.6.1]\label{riganagaga}
The map 
\[X\mapsto X_{\circ}:=\{x\in X | [\HH(x):\K]<\infty\}\]
 is a fully faithful functor from the category of Hausdorff strictly $\K$-analytic spaces to the category of quasiseparated rigid $\K$-spaces. It induces an equivalence between the category of paracompact strictly $\K$-analytic spaces and the category of quasiseparated rigid $\K$-analytic spaces that allow an admissible affinoid cover of finite type.
\end{theorem}

\subsubsection{GAGA results for coherent sheaves and vector bundles}
In this paragraph we review some GAGA results for coherent sheaves and vector bundles. In detail we are concerned with projective schemes and their formal completions (\ref{formelsheafgaga} --- \ref{formelvectgaga}), projective schemes and their rigidifications (\ref{kpflemma} and \ref{kpflemmarem}) and analytic spaces and their associated rigid spaces (\ref{riganavect}). In the last part we introduce the notion of a {\em model} of a vector bundle and prove that a vector bundle has an algebraic model if and only if it has a formal rigid model (\ref{modelgaga}).
\begin{theorem}[\cite{EGA} III Corollaire 5.1.6]\label{formelsheafgaga}
Let $A$ be a Noetherian $I$-adic ring, $\X$ a projective $A$-scheme and $\widehat{\X}$ its formal completion along $\X/I\Loc_{\X}$. Then the following holds: a sheaf on $\widehat{\X}$ is coherent if and only if it is the formal completion of a coherent sheaf on $\X$.
\end{theorem}
\begin{corollar}
Let $\K$ be a local number field, $\K^{\circ}$ its ring of integers and $\X$ a projective $\K^{\circ}$-scheme. Let $\L$ be a complete subfield of $\Cp$ that is an extension of $\K$. Denote its ring of integers by $\L^{\circ}$. Define $\X_{\L^{\circ}}:=\X\otimes\L^{\circ}$ and denote by $\widehat{\X_{\L^{\circ}}}$ its formal completion along $\X_{\L^{\circ}}/p\Loc_{X_{\L^{\circ}}}$. Then the following holds: A sheaf on $\widehat{\X_{\L^{\circ}}}$ is coherent if and only if it is the formal completion of a coherent sheaf on $\X_{\L^{\circ}}$.
\end{corollar}
\begin{proof}
Note that $\widehat{\X_{A}}=\widehat{\X\times\Spec A}=\widehat{\X}\times\Spf A=\widehat{\X}_{A}$ for every continuous morphism $\K^{\circ}\ra A$ of adic rings. 
Let $\E$ be a coherent sheaf on $\widehat{\X_{\L^{\circ}}}$. There are finitely many elements $w_{1},...,w_{n}\in\L^{\circ}$ such that $\E$ comes from a coherent sheaf $\E'$ on $\widehat{\X}_{A}$ (which is equal to $\widehat{\X_{A}}$), where $A$ is the ring of restricted power series $\K^{\circ}\{w_{1},\ldots,w_{n}\}$. 
Define $\X':=\X_{A}$. We are going to show that in the diagram
\[\xymatrix{\X'_{\L^{\circ}}\ar[r]^{f}\ar[d]^{\widehat{\cdot}}&\X'\ar[d]^{\widehat{\cdot}}\\
\widehat{\X'_{\L^{\circ}}} \ar@{-->}[r]^{\widehat{f}}&\widehat {\X'}}\] 
the morphism $\widehat{f}$ exists and that $\widehat{f}^{*}\widehat{\FF}=\widehat{f^{*}\FF}$ is true. The claim follows by Paragraph \cite[10.9.1]{EGAn} if one notes that:
\begin{enumerate}
\item $\X'$ is an $A$-scheme of topological finite presentation, and $\X_{\L^{\circ}}$ is an $\L^{\circ}$-scheme of topological finite presentation.
\item Therefore the defining ideals of the special fibres $p\Loc_{\X'}$ and $p\Loc_{\X'_{\L^{\circ}}}$ are coherent, since they are of finite presentation.
\item It is true that $f^{*}(p\Loc_{\X'})\Loc_{\X'_{\L^{\circ}}}=p\Loc_{\X'_{\L^{\circ}}}$. 
\end{enumerate}
Thus $\E=\widehat{f}^{*}\E'$. By \cite[0 Proposition 7.5.4 (iii)]{EGA} $A$ is adic and Noetherian. The ring $A$ is a subring of $\L^{\circ}$ and the embedding is continuous. Hence by theorem \ref{formelsheafgaga} there is a coherent sheaf $\FF$ on $\X_{A}$ whose formal completion is $\E'$. 
Therefore there is a coherent sheaf $f^{*}\FF$ that satisfies $\widehat{f^{*}\FF}=\E$ 
\end{proof}
\begin{corollar}\label{formelvectgaga}
In the situation of the above theorem and the corollary: a sheaf $\FF$ on $\widehat{\X}$ is locally free if and only if it is the formal completion of a locally free sheaf $\E$ on $\X$.
\end{corollar}
\begin{proof}
By \cite[Proposition 10.10.2.9]{EGAn} locally on $\widehat{\X}$ we have \[\FF_{|\Spf(A)}=\widehat{\widetilde{M}}=\widehat{\E}_{|\Spf(A)}=\widehat{\widetilde{N}}\] if $\E_{|\Spec (A)}=\widetilde{N}$. Since $\widehat{\widetilde{\cdot}}$ is an equivalence of categories between $A$-modules of finite type and coherent $\Loc_{\Spf(A)}$-modules by \cite[I Th\'eor\`eme 10.10.2]{EGAn}, the isomorphism $M\cong N$ follows. The sheaf $\E$ is locally free if and only if $N$ is a projective $A$-module according to \cite[I Corollaire 1.4.4]{EGAn}. But by \cite[I Proposition 10.10.8.6]{EGAn} this is equivalent to $\FF$ being locally free. 
\end{proof}
Now we come to sheaves on rigid spaces. 
\begin{theorem}[\cite{Kpf1974}, 3. Gaga-Satz 5.1]\label{kpflemma}
Let $X$ be a projective variety over $\K$. Then the category of coherent algebraic sheaves on $X$ is equivalent to the category of rigid coherent sheaves on $X^{rig}$.
\end{theorem}
\begin{corollar}\label{kpflemmarem}
The category of algebraic vector bundles on $X$ is equivalent to the category of rigid vector bundles on $X^{rig}$.
\end{corollar}
\begin{proof}
Let $E$ be a coherent sheaf on an algebraic $\K$-variety $X$. Locally on an affine subscheme $U$ of $X$ the rigidification of $E$ is the coherent rigid sheaf associated to the module $E(U)$, that is for every $\Sp(B)$ admissible open in $U^{rig}$ it is 
\[E^{rig}(\Sp B)=E(U)\otimes_{\Loc_{X}(U)}B\]
as proved in \cite[Bemerkung 3.2]{Kpf1974}. $E$ is locally free if and only if for every point in $X$ there is an affine neighbourhood $U$ such that $E(U)$ is a free $\Loc_{X}(U)$-module. This is equivalent to the assertion $E(U)\otimes B$ is a free $B$-module.     
\end{proof}
\begin{proposition}[\cite{Berk1993}, Page 37, Proposition 1.3.4]\label{riganavect}
With the notations in the theorem \ref{riganagaga} the following statement is true: if $X$ is a good strictly $\K$-analytic space, then the categories of coherent sheaves resp. the categories of vector bundles on $X$ and on $X_{\circ}$ are equivalent.
\end{proposition}
\begin{nota}
Let $E$ be a vector bundle on a projective $\K$-curve. Its rank and its degree are equal to the rank and degree of the rigidification of $E$ and of the analytification of $E$. The invariance of the degree is due to the fact that the cohomology group $H^{n}(X,E)$ is isomorphic to $H^{n}(X^{rig},E^{rig})$ and $H^{n}(X^{an},E^{an})$. In particular $E$ is semistable if and only if its rigidification resp. its analytification is semistable.
\end{nota}

Let $\K$ be a complete subfield of $\Cp$ and $\K^{\circ}$ its ring of integers. Let $X$ be a projective and smooth curve over $\K$, $\X$ a $\K^{\circ}$-model of $X$ with the embedding $i: X \inj \X$ and $\widehat{\X}$ its formal completion along $\X/p\Loc_{\X}$. For a vector bundle $\E$ on $\X$ the term $\widehat{\E}$ denotes the formal completion of $\E$ along the closed fibre defined be the ideal $p\Loc_{\X}$.
\begin{defin}
Let E be a vector bundle on $X$. If there is a vector bundle $\E$ on $\X$ satisfying $E\cong i^{*}\E$, then $\E$ is called ({\em algebraic $\X$-}){\em model} of $E$.\\
Let $(X,r,\widetilde{X})$ be a formal rigid space. If there is a formal rigid vector bundle $\E$ on $(X,r,\widetilde{X})$, such that \[E(r^{-1}(U))\cong\E(U)\otimes\K\] for every open affine $U\subset \overline{X}$, then $\E$ is called {\em formal} ({\em $(X,r,\widetilde{X})$-}){\em model} of $E$.
\end{defin}
\begin{lemma}\label{modelgaga}
Let $E$ be a vector bundle on $X$ and $\E'$ a vector bundle on $\X$. The following statements are equivalent:
\begin{enumerate}
\item $E\cong i^{*}\E'$
\item $r_{*}E^{rig}\cong\widehat{\E'}\otimes\K$.
\end{enumerate}
\end{lemma}
\begin{proof}
To prove $(1)\Ra (2)$ choose a covering of $\X$ by open affines that trivialises $\E$. Assume for an open affine $U=\Spec A$ of the covering that we have $\E'_{|U}=\oplus_{j=1}^{n}f_{j}\Loc_{U}$. It follows
\begin{eqnarray*}
\widehat{\E'_{|U}}\otimes\K
&=&\widehat{\oplus_{j=1}^{n}f_{j}\Loc_{U}}\otimes\K
=\oplus_{j=1}^{n}f_{j}\widehat{\Loc_{U}}\otimes\K\\
&=&\oplus_{j=1}^{n}f_{j}r_{*}\Loc_{\Spec(A\otimes\K)}^{rig}
=r_{*}(\oplus_{j=1}^{n}f_{j}i^{*}\Loc_{U})^{rig}=r_{*}(i^{*}\E'_{|U})^{rig}\\
&\cong&r_{*}E_{|U}^{rig}.
\end{eqnarray*}
The last isomorphism comes from $E\cong i^{*}\E'$, hence the local isomorphisms can be glued to a global isomorphism.\\
To prove $(2)\Ra (1)$ we define $E':=i^{*}\E'$ and note that by assumption \[r_{*}E^{rig}\cong\widehat{\E'}\otimes\K\cong r_{*}E'^{rig}.\] Since we only have to prove that for any two vector bundles $E_{1}$ and $E_{2}$ that satisfy $r_{*}E_{1}^{rig}=r_{*}E_{2}^{rig}$, the assertion $E_{1}\cong E_{2}$ is true. This follows by \ref{redaffine} in definiton \ref{reddef} and theorem \ref{kpflemma}.
\end{proof}

\subsubsection{Quotients by finite groups}
\begin{nota}
Let $G$ be a finite group acting on a separated scheme $X$. Suppose that every point $x\in X$ has an affine neighbourhood that is stable under $G$. Then the quotient $X/G$ in the category of schemes exists. For instance this is the case if $X$ is quasi-projective. In the following, when we write {\em the quotient $X/G$ exists} we mean that {\em every point $x\in X$ has an affine neighbourhood that is stable under $G$}.
\end{nota}

\begin{nota}
As stated in \cite[Ex 4.3.18]{Liu2002}, schematic quotients commute with flat base change. 
\end{nota}
\begin{lemma}
Let $A$ be a Noetherian commutative ring that is flat over a discrete valuation ring $R$ with uniformising parameter $\pi$, and let $I:=\pi A$ be the ideal of $A$ that is generated by $\pi$. Denote by $\widehat{A}$ the completion of $A$ in the $I$-adic topology. If a finite group $G$ acts $R$-linearly on $A$, then by functoriality of $ \widehat{\cdot}$ it also acts on $\widehat{A}$, and the fixed rings satisfy 
\[\widehat{A}^{G}=\widehat{A^{G}},\] 
where $\widehat{A^{G}}$ denotes the completion in the $(I\cap A^{G})$-adic topology.
\end{lemma}
\begin{proof}
First we have to prove that on $A^{G}$ the induced topology and the $(I\cap A^{G})$-adic topology coincide. 
It is true that $I^{n}=\pi^{n} A$ and we have $(\pi A)\cap A^{G}=\pi A^{G}$, because $A$ has no $\pi$-torsion, since it is $R$-flat. The equality $((\pi A)\cap A^{G})^{n}=(\pi A^{G})^{n}=\pi^{n}A^{G}$ follows. Hence both topologies coincide.\\
In the category of topological rings it is true that $A^{G}\subset A$ and $\widehat{A^{G}}\subset \widehat{A}^{G}$ as well, since the morphism of topological rings 
\[\widehat{\phi}:\widehat{A}\ra\prod_{g\in G} \widehat{A}\text{, } a\mapsto(ga-a)_{g\in G}\]
is continuous. Therefore it is enough to show that $\widehat{A^{G}}$ and $\widehat{A}^{G}$ coincide as abelian groups.\\
Let $\phi$ be the morphism $\phi(a):=(ga-a)_{g\in G}$ and consider the exact sequence 
\[\xymatrix{0\ar[r]& A^{G}\ar[r]& A\ar[r]^{\phi}&{\prod_{g\in G} A}}\]
of abelian groups, where the group $A$ in the middle induces the topologies of the right and left term. By \cite[Corollary 10.3]{Atiy1969} formal completion is an exact functor; it follows that the sequence
\[\xymatrix{0\ar[r]& {\widehat{A^{G}}}\ar[r]& {\widehat{A}}\ar[r]^{\widehat{\phi}}&{\prod_{g\in G} \widehat{A}}}\] is exact.
Therefore $\widehat{A^{G}}$ is the kernel of $\widehat{\phi}$. By what we have proved above, on the left term the induced topology and the $I\cap A^{G}$-adic topology coincide. The claim follows.
\end{proof}
\begin{lemma}\label{formalquot}
Let $G$ be a finite group and $R$ a discrete valuation ring with uniformising element $\pi$. Let $Y$ be a flat, separated and Noetherian $R$-scheme endowed with an $R$-linear action of $G$, such that the schematic quotient $X:=Y/G$ exists, and denote formal completion along the closed fibre by $\widehat{\cdot}$. Then $G$ acts on $\widehat{Y}$ by functoriality and it is true that \[\widehat{X}=\widehat{Y}/G\] holds.
\end{lemma}
\begin{proof}
Assume $Y$ to be affine. Then there is a ring $A$ that is flat over $R$, such that $Y=\Spec(A)$ and $X=\Spec(A^{G})$. Denote by $Y_{s}$ the special fibre of $Y$. If we define $I:=\pi A$, then $Y_{s}=\Spec(A/I)=\Spec(A/\pi A)$. It follows 
\[X_{s}=\Spec(A^{G}/(I\cap A^{G})).\] To prove that $\widehat{Y}/G$ and $\widehat{Y/G}$ coincide, we have to show that $\widehat{A}^{G}$ (completed in the $I$-adic topology) coincides with $\widehat{A^{G}}$ (completed in the $I\cap A^{G}$-adic topology). This was proved in the lemma above.\\
For the general case choose an affine cover of $Y$ by $G$-invariant affine subschemes such that all intersections of the covering sets are also affine (such a covering exists because $Y$ is separated and the quotient exists). Then the formal quotient can be glued from the local pieces. Since $\widehat{X}$ and $\widehat{Y}/G$ coincide locally, they coincide globally.
\end{proof}
\begin{lemma}
Let $R$ be a discrete valuation ring and $\Y$ a separated and admissible formal $R$-scheme endowed with an $R$-linear action by a finite group $G$ such that the quotient $\X=\Y/G$ of formal $R$-schemes exists. Then $G$ acts on $\Y^{rig}$ by functoriality and \[\X^{\text{rig}}=\Y^{\text{rig}}/G.\]
\end{lemma}
\begin{proof}
If $\Y$ is affine, then there is an admissible $R$-algebra $A$ such that $\Y=\Spf(A)$ and $\X=\Spf(A^{G})$. Let $Q$ be the quotient field of $R$. By application of the functor $\cdot^{\text{rig}}$ from formal to rigid spaces we obtain \[\Y^{rig}=\Spm(A\otimes Q) \text{ and }\X^{rig}=\Spm(A^{G}\otimes Q).\]
By \cite[0 Proposition 7.6.13]{EGA} the morphism $A\ra A\otimes Q$ is flat; taking fixed points commutes with flat base change. Thus $A^{G}\otimes Q=(A\otimes Q)^{G}$,  and we have 
\[\X^{\text{rig}}=\Spm((A\otimes Q)^{G}).\] Therefore $\X^{\text{rig}}=\Y^{\text{rig}}/G$ in the affine case.\\
In the general case there is a $G$-invariant affine covering of $\Y$ that has affine intersections. Then one glues the local pieces and obtains a global quotient.
\end{proof}
\subsection{Fundamental groups of analytic manifolds}
In this subsection we quickly review the general theory of fundamental groups. We introduce some fundamental groups of importance for our work and gather their relations. We define the fundamental group of finite topological coverings, which is best suited for our purpose. We prove that this group satisfies Grothendieck's six axioms of a Galois theory.
The presentation refers to Yves Andr\'e's book \cite{Andr2003a}.\\ 
In the following let $\K$ be a subfield of $\Cp$ that is complete with respect to the induced norm. 
\subsubsection{Fundamental groups}\text{ }\\
We give a brief review of the theory of fibre functors and fundamental groupoids and define some interesting fundamental groups. 
\begin{defin}[\cite{Andr2003b}, Paragraph 4.1]
A smooth paracompact strictly $\K$-analytic space $S$ is called {\em analytic $\K$-manifold}.
\end{defin}
\begin{nota}
Note that every analytic $\K$-manifold is a good analytic space, since it is smooth over $\M(\K)$, which is a good analytic space.
\end{nota}
\begin{theorem}[\cite{Berk1999}, Theorem 9.1]
Analytic manifolds are locally contractible. Therefore every analytic manifold has a universal covering.
\end{theorem}

\begin{defin}\text{ }
\begin{enumerate}
\item A morphism of analytic manifolds $f:S'\ra S$ is called a {\em covering} of $S$ if 
\begin{enumerate}
\item $S$ is covered by a family of open sets $(U_{i})_{i\in I}=:\U$, 
\item for all sets $U$ in $\U$ it is true that \[f^{-1}(U)=\coprod_{j\in J}V_{j},\] 
\item and for all $j\in J$ the restriction $f_{|V_{j}}:V_{j}\ra U$ is finite. 
\end{enumerate}
The family $\U$ will be called {\em covering collection}. A morphism of coverings $S'\ra S$ to $S''\ra S$ is a commutative diagram
\[
\xymatrix{S' \ar[rr] \ar[rd] && S'' \ar[ld]\\
&S\text{ .}}\]
The thereby defined {\em category of coverings of $S$} is called $\Cov_{S}$.\\
\item By $\Cov^{et}_{S}$ we denote the {\em \'etale} coverings: a covering $f:S'\ra S$ is \'etale if there is a covering collection $\U$ of $S$ such that the restrictions $f_{|V_{j}}:V_{j}\ra U$ are \'etale for all $U\in\U$.\\
\end{enumerate}
\end{defin}
\begin{defin}
Let $S'\ra S$ be an \'etale covering and $R\subset S'\times_{S} S'$ an equivalence relation on $S'$ over $S$, which is a union of connected components. The quotient $S'/R$, viewed as an \'etale sheaf as in \cite[Expos\'e V 1 ]{SGA1}, is representable by an \'etale covering $S''\ra S$ as proved in \cite[Lemma 2.4]{Jong1995}. It is called {\em \'etale quotient covering} or just {\em $S$-quotient}.
\end{defin}

\begin{defin}\label{covtypes}
There are the following interesting full subcategories of $\Cov_{S}^{et}$, which are stable under connected components, $S$-fibre products and $S$-quotients.
\begin{enumerate}
\item $\Cov^{alg}_{S}$: a finite morphism $f:S'\ra S$ that is an \'etale covering is called {\em finite \'etale} covering. Finite \'etale coverings are also called {\em algebraic coverings} because of the Riemann existence theorem (see below and theorem \ref{curveana}).\\
\item $\Cov^{top}_{S}$: a covering $f:S'\ra S$ is {\em topologic} if there is a covering collection $\U$ of S such that the restrictions $f_{|V_{j}}:V_{j}\ra U$ are isomorphisms for all $U$ in $\U$.\\
\item $\Cov^{ftop}_{S}$: a finite morphism $f:S'\ra S$ that is a topological covering is called {\em finite topological} covering.\\
\item $\Cov^{temp}_{S}$: a covering $f:S'\ra S$ is called {\em temperated} if it is an $S$-quotient of a composite \'etale covering $T'\ra T\ra S$ where $T'\ra T$ is a topological covering and $T\ra S$ an algebraic covering.
\end{enumerate}
We use $\Cov_{S}^{\bullet}$ if we mean any of these categories.
\end{defin}

\begin{defin}
A geometric point of $S$ is a morphism $\overline{s}: \M(\Cp)\ra S$. If $\overline{s}$ is a geomteric point, then we denote the unique point in its image in $S$ by $s$.
\end{defin}
\begin{defin}\text{ }
\begin{enumerate}
\item Let $\overline{s}$ be a geometric point of $S$. The set 
\[f^{-1}(\overline{s}):=\{\overline{s'}:\M(\Cp)\ra S'|f\circ \overline{s'}=\overline{s}\}\] 
is called {\em geometric fibre} of $\overline{s}$ in the covering $f:S'\ra S$.
The covariant functor 
\[F_{S,\overline{s}}^{\bullet}:\Cov_{S}^{\bullet}\ra\Ens, (f:S'\ra S)\mapsto f^{-1}(\overline{s})\] is called {\em fibre functor}.\\
\item A {\em $\bullet$-path} (\'etale path, finite \'etale path, and so on) from a geometric point $\overline{s}$ to a geometric point $\overline{t}$ is an isomorphism of fibre functors $F_{S,\overline{s}}^{\bullet}\cong F_{S,\overline{t}}^{\bullet}$.
The set of $\bullet$-paths can by topologised by taking as the fundamental open neighbourhoods of a path $\alpha$ the sets $\Stab_{S',\overline{s'}}\circ\alpha$ where $\Stab_{S',\overline{s'}}$ runs among the stabilisers in $\Aut(F^{\bullet}_{S,\overline{s}})$ of arbitrary geometric points $\overline{s'}$ above $\overline{s}$ in arbitrary $\bullet$-coverings $S'\ra S$ in $\Cov^{\bullet}_{S}$.\\
\item If $f:S\ra T$ is a morphism of $\K$-manifolds and $\alpha$ a $\bullet$-path from $\overline{s}$ to $\overline{t}$, then $f\circ\alpha$ is a $\bullet$-path from $f\circ \overline{s}$ to $f\circ \overline{t}$. This is compatible with composition of paths.
\end{enumerate}
\end{defin}

Comparison to algebraic theory: let $X$ be a smooth $\K$-Variety, its analytification $X^{an}$ is a $\K$-manifold by \cite[III Remark 1.1.2 e)]{Andr2003a}. Via the analytification functor the category of algebraic coverings -- which has been defined by Grothendieck in \cite{SGA1} -- is equivalent to the category $\Cov_{X^{an}}^{alg}$. This is a consequence of the Riemann existence theorem proved in \cite{Ltke1993}. Therefore algebraic paths exist on $X^{an}$. Due to de Jong (\cite[Theorem 2.0]{Jong1995}) \'etale paths exist as well. Every topological path between $s$ and $t$ lifts to an \'etale paths between $\overline{s}$ and $\overline{t}$.
\begin{defin}\label{funddef}\text{ }
\begin{enumerate}
\item The category whose objects are the geometric points of $S$ and whose homomorphisms are \[\Hom^{\bullet}(\overline{s},\overline{t}):=\Iso(F_{S,\overline{s}}^{\bullet},F_{S,\overline{t}}^{\bullet})\] 
is called the {\em $\bullet$-fundamental groupoid} and is denoted by $\Pi_{1}^{\bullet}(S)$.
\item The group $\pi_{1}^{\bullet}(S,\overline{s}):=\Aut(F^{\bullet}_{S,\overline{s}})$ is called {\em $\bullet$-fundamental group of $S$ with basepoint $\overline{s}$}.
\end{enumerate}
\end{defin}
\begin{nota}[\cite{Jong1995}, Lemma 2.6]
The category of topological coverings of $S$ is equivalent to the category of covering spaces of the underlying topological space $|S|$ of $S$. The topological fundamental group $\pi_{1}^{top}(S,\overline{s})$ does only depend on $|S|$ and on the unique point $s$ in the image of the geometric point $\overline{s}$ in $S$.
\end{nota}

\subsubsection{The finite topological fundamental group}\text{ }\\
In \cite[V.4]{SGA1} Grothendieck defined six axioms for a Galois theory; we state them and prove that for a connected base $S$ the category of finite topological coverings $\Cov^{ftop}_{S}$ together with the fibre functor $F^{ftop}$ satisfies these axioms.
\begin{defin}\label{galoisaxioms}
Let $\C$ be a category and $F$ a covariant functor \[F:\C\ra\Ens.\] The category $\C$ together with the functor $F$ is called a {\em Galois theory} if the following six axioms are satisfied:
\begin{enumerate}
\item[(G1)] The category $\C$ has a final object, and the fibre product of two objects over a third object exists in $\C$.
\item[(G2)] Finite sums exist in $\C$. The quotient of an object by a finite group of automorphisms exists.
\item[(G3)] Let $u:S'\ra S''$ be a morphism in $\C$, then $u$ factorises in a composition 
\[\xymatrix{S'\ar[r]^{u'}& S''_{1}\ar[r]^{u''}& S''},\]
where $u'$ is a strict epimorphism and $u''$ is a monomorphism. Furthermore $u''$ is an isomorphism onto a direct summand of $S''$.
\item[(G4)] The functor $F$ is left exact, that is $F$ maps monomorphisms to monomorphisms and commutes with fibre products.
\item[(G5)] The functor $F$ commutes with finite direct sums, maps strict epimorphisms to epimorphisms and commutes with quotients by finite groups of automorphisms.
\item[(G6)] Assume $u:S'\ra S''$ is a morphism in $\C$ such that $F(u)$ is an isomorphism, then $u$ itself is an isomorphism.
\end{enumerate}
\end{defin}
In order to prove the property (G1) for $F^{ftop}$ we need the following lemma.
\begin{lemma}\label{topcovering}
Let $S$ be a connected $\K$-manifold and $f:S'\ra S$ and $g:S''\ra S$ be two topological coverings. Every $S$-morphism $\phi:S'\ra S''$ is a topological covering. 
\end{lemma}
\begin{proof}
Let $\pi:\widetilde{S}\ra S$ be the universal covering space of $S$. Thus there is a covering $(U_{i})_{i\in I}$ of $S$ such that for every $U_{i}$ exists a family $(V_{ij})_{j\in J}$ in $\widetilde{S}$ such that
\[\pi^{-1}(U_{i})=\coprod_{j} V_{ij}\text{, and }\pi_{|V_{ij}}:V_{ij}\ra U_{i}
\;\text {is an isomorphism.}\] Since $\widetilde{S}$ is the universal covering space of $S$, there are morphisms of analytic spaces 
\[\xymatrix{{\pi:\widetilde{S}\ar[r]^{\pi'}} & S' \ar[r]& S},\] 
where $S'=\widetilde{S}/H$ and $H$ is the fundamental group of $S'$. Thus there is a decomposition
\[\xymatrix{{\pi_{|\pi^{-1}(U_{i})}:\coprod_{j} V_{ij}\ar[r]}& (\coprod V_{ij})/H \ar[rr]^{f_{|f^{-1}(U_{i})}} & & U_{i}}.\] 
Up to isomorphism the quotient $\coprod V_{ij}/H$ identifies some of the $V_{ij}$. We denote the images of $V_{ij}$ in the quotient $\widetilde{S}/H$ by $V'_{ij}$. Hence we have the decomposition
\[\xymatrix{\pi_{|\pi^{-1}(U_{i})}:\coprod_{j} V_{ij}\ar[r]& \coprod_{k} V'_{ij_{k}} \ar[rr]^{f_{|f^{-1}(U_{i})}}& & U_{i}},\] 
where $f_{|V'_{ij_{k}}}:V'_{ij_{k}} \ra U_{i}$ are isomorphisms.\\
With the analogous labeling the analogous statement is true for $S''$, that is the morphisms $g_{|V''_{ij_{m}}}: V''_{ij_{m}}\ra U_{i}$ are isomorphisms.\\
Therefore we have the following diagram
\[
\xymatrix{V'_{ij} \ar[rr]^{\phi} \ar[rd]_{f_{|V'_{ij}}} && V''_{ik} \ar[ld]^{g_{|V''_{ik}}}\\
&U_{i}.}\]
The morphism $\phi$ is an isomorphism, because $f_{|V'_{ij}}$ and $g_{|V''_{ik}}$ are isomorphisms. Since $S''$ is covered by the sets $V''_{ik}$, the morphism $\phi$ is a topological covering. 
\end{proof}
\begin{proposition}
Let $S$ be a connected and proper $\K$-manifold the category of finite topological coverings $\Cov^{ftop}_{S}$ together with the fibre functor $F^{ftop}$ satisfies the axioms of a Galois theory as defined in \ref{galoisaxioms}.
\end{proposition}

\begin{proof}\text{ }
 \begin{enumerate}
\item[(G1)] The category $\Cov^{ftop}_{S}$ admits $S$-fibre products. By the lemma above a general fibre product $S'\times_{S'''}S''$ in $\Cov^{ftop}_{S}$ is just the $S'''$-fibre product in $\Cov^{ftop}_{S'''}$. The final object exists.
\item[(G2)] This is true for topological coverings. Finite coverings are closed under finite sums and quotients by finite groups. Therefore the assertion is true for finite topological coverings. 
\item[(G3)] To prove this property we use the notations of Lemma \ref{topcovering}. Since $f:S'\ra S$ and $g:S''\ra S$ are \'etale coverings, they are \'etale morphisms. The reason for this is that there is a covering $(U_{i})_{i\in I}$ of $S$ such that for all $U_{i}$ there is a covering $(V'_{ij})_{j\in J}$ of $S'$ with the property
\[f^{-1}(U_{i})=\coprod_{j} V'_{ij}\] 
and \[f_{|V'_{ij}}:V'_{ij}\ra U_{i}\]
is \'etale.
It follows that 
\[f_{|f^{-1}(U_{i})}:\coprod_{j} V'_{ij}\ra U_{i}\] 
is \'etale. \'Etaleness is a local condition, therefore the morphism $f:S'\ra S$ is \'etale. For the morphism $g$ it is the same proof.\\
Let $u:S'\ra S''$ be a morphism in $\Cov^{ftop}_{S}$. We use the equivalence between categories of finite \'etale (analytic) coverings and algebraic coverings (\ref{curveana}) and conclude, that the morphism $f^{alg}$ is \'etale. Hence the image $u^{alg}(S'^{alg})$ is open and closed in $S''^{alg}$; therefore it is a connected component of $S''^{alg}$, and $f^{alg}$ factorises in a strict epimorphism $u'^{alg}$ and a monomorphism $u''^{alg}$ as proved in \cite[V Proposition 3.5]{SGA1}. By equivalence of categories this also holds for $f$ itself. Hence we get the following diagram.  
\[\xymatrix{S' \ar@/^1.5pc/[rr]^{u}\ar[r]^{u'} \ar[rd]_{f} & f(S') \ar[r]^{u''}\ar[d]^{h} & S'' \ar[ld]^{g}\\ & S .&}\]
The space $S''$ has a decomposition $f(S')\coprod S_{1}''$, where $f(S')$ is a connected component of $S''$.
Notice here that an equivalence of categories maps monomorphisms to monomorphisms, epimorphisms to epimorphisms and strict epimorphisms to strict epimorphisms. This is due to the fact that these notions are defined by universal properties.\\
The morphism $h$ is finite. We have to show that it is a topological morphism. As in the proof of the above lemma we find a covering $(U_{i})_{i\in I}$ of $S$ such that 
\[f^{-1}(U_{i})= \coprod_{j}V'_{ij}\text{ and }\] \[g^{-1}(U_{i})=\coprod_{k}V''_{ik}.\] 
Thus we have the equation
\[
\begin{array}{ll}\nonumber
g^{-1}(U_{i})&\cong\coprod_{k}V_{ik}\cong(f(S')\coprod S''_{1})\cap g^{-1}(U_{i})\\
\nonumber &= (f(S')\cap g^{-1}(U_{i}))\coprod (S''_{1}\cap g^{-1}(U_{i})). 
\end{array}
\]
Here remember that $V'_{ij}\cong V_{ij}\cong V''_{ij}$.
The $V_{ij}$ are connected  (since $S$ is connected) therefore 
\[h^{-1}(U_{i})=f(S')\cap g^{-1}(U_{i})\cong\coprod_{k} V_{ij_{k}}.\]Thus $h:f(S')\ra S$ is a finite topological covering. 
\item[(G4)] The algebraic coverings satisfy this axiom, there is an equivalence of categories between finite \'etale and algebraic coverings, a bijection between geometric points, and the fibre functor is defined in the same way. Hence this axiom is also satisfied for the (finite \'etale and therefore for the) finite topological coverings.
\item[(G5)] The same reason as for G4.
\item[(G6)] The same reason as for G4.
\end{enumerate}
\end{proof}
\begin{corollar}
Every finite topological covering is dominated by a (finite) Galois covering and $\pi^{ftop}_{1}(S,\overline{s})$ is a pro-finite group.
\end{corollar}
\begin{proof}
This follows by the general Galois theory of Grothendieck (cf. \cite[page 123]{SGA1}). 
\end{proof}

\begin{nota}\label{fundamentalfactor}\text{ }
\begin{enumerate}
\item The fundamentel group $\pi^{ftop}(S,\overline{s})$ only depends on the underlying topological space $|S|$ of $S$ and on the unique point $s$ in the image of $\overline{s}$.  
\item The fully faithful embedding functor $H:\Cov^{ftop}\inj \Cov^{alg}$ satisfies $F^{ftop}=F^{alg}\circ H$. Because $H$ is fully faithful, the morphism 
\[\pi_{1}^{alg}(S,\overline{s})\ra\pi_{1}^{ftop}(S,\overline{s})\]
 is surjective by \cite[V Proposition 6.9]{SGA1}.
\item Furthermore $F^{alg}_{X_{\Cp},\overline{x}}$ and $F^{ftop}_{X_{\Cp},\overline{x}}$ are strictly pro-representable by $(Y_{i},y_{i},\phi_{ij})$ and $(Y'_{i},y'_{i},\phi'_{ij})$ respectively. If $Z\ra X$ is a finite topological covering and $\overline{z}\in Z(\Cp)$ is a geometric point, then the following diagram is commutative:
\[\xymatrix{{\pi^{alg}(X_{\Cp},\overline{x})\ar[r]^{\phi_{y_{i}}}\ar@{->>}[d]}& {\Gal_{X_{\Cp}}Y_{i}\ar[r]^{\phi_{\overline{z}}}\ar[d]^{\phi_{y'_{i}}}}&{\Gal_{X_{\Cp}}Z_{\Cp}}\\
{\pi^{ftop}(X_{\Cp},\overline{x})\ar[r]^{\phi_{y'_{i}}}}&{\Gal_{X_{\Cp}}Y'_{i}\ar[ru]_{\phi_{\overline{z}}}} &}.\]
\item Because of the Riemann existence theorem, $\Pi_{1}^{alg}(S)$ and $\pi_{1}^{alg}(S,\overline{s})$ respectively agree with their counterparts of \cite{SGA1} and therefore with those in \cite{Deni2004}.
\item The inclusion $\Cov^{top}_{S}\inj\Cov^{temp}_{S}$ gives rise to a surjective homomorphism \[\pi_{1}^{temp}(S,\overline{s})\proj\pi_{1}^{top}(S,\overline{s}).\]
The surjectivity follows by \cite[III Corollary 1.4.8]{Andr2003a} (that is $\pi_{1}^{temp}(S,\overline{s})\ra\pi_{1}^{top}(S,\overline{s})$ has dense image) and the fact that $\pi_{1}^{top}(S,\overline{s})$ is discrete.
\item By a remark by Yves Andr\'e in \cite[Page 128]{Andr2003a} the groups $\pi_{1}^{alg}(S,\overline{s})$ and $\widehat{\pi^{temp}(S,\overline{s})}$ can be identified.
\item Since $\Cov_{S}^{ftop}$ is a full subcategory of $\Cov_{S}^{top}$ and $\pi^{ftop}_{1}(S,\overline{s})$ is pro-finite, we have by \cite[III Corollary 1.4.8]{Andr2003a} that 
\[\pi_{1}^{ftop}(S,\overline{s})=\widehat{\pi_{1}^{top}(S,\overline{s})}.\]
\end{enumerate}
\end{nota}
\subsection{Mumford curves}\label{mumford.section}
In this subsection we review the definition and some facts about Mumford curves and present two examples.
\subsubsection{Definition and facts}\text{ }\\ 
Let $\K$ be a complete subfield of $\Cp$.
\begin{defin} Let $\PP$ be the rigid projective line over $\K$, and let $\Gamma$ be a subgroup of $\PGL_{2}(\K)$.
The group $\Gamma$ acts on $\PP$.
\begin{enumerate}
\item A point $x\in\PP$ is called a {\em limit point} of $\Gamma$ if there is a sequence
$(\gamma_{n}\in\Gamma)_{n\in\N}$ of pairwise different elements of $\Gamma$ and a point $y\in\PP$ such that $\lim_{n\ra\infty}\gamma_{n}(y)=x$. The set of all limit points of $\Gamma$ will be
called $\Sigma_{\Gamma}$.
\item The group $\Gamma$ is called {\em discontinuous} if $\Sigma_{\Gamma}$ is not equal to $\PP$ and if for every point $x\in\PP$
the closure of the orbit $\Gamma x$ is compact. If $\K$ is locally compact, then the last condition is trivially satisfied.
\item The group $\Gamma$ is called {\em Schottky group} if it is finitely generated, discontinuous and has no elements of finite order different from 1.
\end{enumerate}
\end{defin}
\begin{defin} Let $\Gamma$ be a Schottky group.  Then $\Gamma$ acts freely on the open set
$\Omega_{\Gamma}:=\PP-\Sigma_{\Gamma}$.  The quotient $\Omega_{\Gamma}/\Gamma$ is compact.  Endowed with the obvious $\K$-rigid
structure it becomes a rigid $\K$-curve. By \cite[III Theorem 2.2]{Gerr1980} it is the rigidification of a projective and smooth $\K$-curve $X_{\Gamma}$.  The curve $X_{\Gamma}$ is called the {\em Mumford curve} associated to $\Gamma$.
\end{defin}
In this thesis we will solely deal with Mumford curves that can be defined over a local number field. We will also investigate Mumford curves over $\Cp$, but only those that come from a local number field by base extension. For examples of Mumford curves see paragraph \ref{subsexamples}.
\begin{defin}[\cite{Mumf1972}, pages 160 and 164]\text{ }
\begin{enumerate}
\item A {\em stable curve} over a discrete valuation ring $R$ is a proper and flat $R$-scheme whose geometric fibres are reduced, connected and 1-dimensional, have at most ordinary double points, and such that their non-singular rational components, if any, meet the remaining components in at least 3 points.
\item Denote by $k$ the residue field of $R$. A stable curve $\X$ will be called {\em degenerated} if the normalisations of all the components of $\X\otimes_{R}\overline{k}$ are rational curves. Here $\overline{k}$ denotes the algebraic closure of the residue field of $R$.
\item A stable curve $\X$ is called {\em $k$-split degenerated} if the normalisations of all the components of $\X\otimes_{R}k$ are isomorphic to $\PP_{k}^{1}$, and if all the double points are $k$-rational with two $k$-rational branches.
\item A {\em semistable} curve is a flat and proper $R$-scheme whose geometric fibres are reduced, connected and 1-dimensional with at most ordinary double points.
\end{enumerate}
\end{defin} 
\begin{theorem}[\cite{Mumf1972}, Theorem 3.3, Theorem 4.20, Corollary 4.2]
Let $\K$ be a local number field. If $\Gamma$ has $g\geq 1$ generators, then $X_{\Gamma}$ is a smooth curve of genus $g$ and has a stable $\widetilde{\K}$-split degenerated model $\X$ over $\K^{\circ}$, which is uniquely defined up to canonical isomorphism. Vice versa every smooth algebraic $\K$-curve $X$ of genus $g\geq 1$ that has a stable $\widetilde{\K}$-split degenerated model is isomorphic to a Mumford curve $X_{\Gamma}$ for a unique Schottky group $\Gamma$ with $g$ generators.
\end{theorem}
\begin{nota}
In Berkovich theory one has the analogous definition and results. Note here that $\Gamma$ only operates on the points $\PP^{an}(\Cp)$, where $\PP^{an}$ is the analytic projective line. The quotient is the analytification of an algebraic curve by theorem \ref{curveana}.    
\end{nota}
Because Berkovich spaces are topological spaces, the following equivalence makes sense where $\cdot^{an}$ denotes analytification of varieties.
\begin{theorem}[\cite{Berk1990}, Theorem 4.4.1] Let $X$ be a smooth and projective curve over $\K$.  The following properties of $X$ are equivalent:
\begin{enumerate} 
\item The curve $X$ is a Mumford curve.  
\item There is a subset $\Sigma\subset\PP^{an}(\K)$ such that the universal covering $\Omega$ of $X^{an}$ is isomorphic to $\PP^{an}-\Sigma$.
\end{enumerate} 
\end{theorem} 
\begin{corollar}\label{fundextension}
If $X$ is a $\K$-Mumford curve and $\L\subset\Cp$ a complete extension field of $\K$, then \[\pi_{1}^{top}(X^{an},\overline{x})=\pi_{1}^{top}(X^{an}_{\L},\overline{x})\text{ and } \pi_{1}^{ftop}(X^{an},\overline{x})=\pi_{1}^{ftop}(X^{an}_{\L},\overline{x}).\]
\end{corollar}
\begin{proof}
It is enough to prove the assertion for the topological fundamental group. Because $X$ is a Mumford curve we have \[X^{an}=\Omega/\pi_{1}^{top}(X^{an},\overline{x}).\] 
Thus \[X^{an}_{\L}=\Omega_{\L}/\pi_{1}^{top}(X,\overline{x})\] follows by base change. It is enough to prove that $\Omega_{\L}$ is the universal covering space of $X^{an}$. It is a covering space by definiton and it is universal because it is contractible as proved in \cite[Theorem 4.2.1]{Berk1990}.
\end{proof}
The following lemma shows that every finite topological covering of a Mumford curve that is induced by a cofinite normal subgroup of $\Gamma$ is again a Mumford curve. 
\begin{lemma}\label{mumfordcover} 
Given a Mumford curve $\Omega/\Gamma$ of genus $g$ over a local number field $\K$ and a cofinite
normal subgroup $N$ of $\Gamma$, then $\Omega/N$ is a Mumford curve of genus
\[(g-1)\cdot |\Gamma/N|+1,\] and the map
\[\Omega/N\ra\Omega/\Gamma\] is a finite topological covering of degree $|\Gamma/N|$.
\end{lemma} 
\begin{proof} 
As $\Gamma$ is a Schottky group it has no elements of finite order different from 1, it is finitely generated, a free subgroup of $\PGL_{2}(\K)$ and operates discontinuously.  Since $N$ is a subgroup of
$\Gamma$, it also has no elements of finite order different from 1 and operates discontinuously. To prove that $N$ is a Schottky group we only have to show that it is finitely generated.  Since $\Gamma$ is a free group and $N$ is cofinite, it is free and generated by $(g-1)\cdot|\Gamma/N|+1$ elements. This is proved in \cite[\S 1 Example 1, \S 2 Theorem 1]{Ihar1966}. Therefore $N$ is a Schottky group of rank $(g-1)\cdot|\Gamma/N|+1$.\\ 
To prove that
$\Omega/N$ is a Mumford curve, we have to show that $N$ and $\Gamma$ have the same limit points.  Let $q$ be a limit point of $\Gamma$;
in formal terms this is:  \[\exists\; p\in\PP^{1}, (\gamma_{i}\in\Gamma)_{i\in I}:  \lim_{i\ra\infty}\gamma_{i}(p)=q, (\gamma_{i}=\gamma_{j}\Ra i=j).\] Because
of cofiniteness, there is a coset $N\overline{\gamma}$ of $N$ in $\Gamma$ such that
\begin{align*} 
& \exists\; p\in\PP^{1},(\gamma_{i}\in N\overline{\gamma})_{i\in I}:  \lim_{i\ra\infty}\gamma_{i}(p)=q, (\gamma_{i}=\gamma_{j}\Ra i=j)\\ 
\Rightarrow &\exists\; p\in\PP^{1}, (\gamma_{i}\in N\overline{\gamma})_{i\in I}: \lim_{i\ra\infty}(\gamma_{i}\gamma^{-1})(\gamma(p))=q, (\gamma_{i}=\gamma_{j}\Ra i=j)\\
\Rightarrow &\exists\; p'\in\PP^{1}, (\gamma_{i}\in N)_{i\in I}:  \lim_{i\ra\infty}(\gamma_{i})(p')=q, (\gamma_{i}=\gamma_{j}\Ra i=j).
\end{align*}
Thus every limit point of $\Gamma$ is a limit point of $N$.  The other direction is obvious.  Therefore $\Omega/N$ is a Mumford curve.\\
Its genus is equal to the number of generators of $N$, hence it is equal to $(g-1)\cdot |\Gamma/N|+1$. 
In Berkovich's
sense the morphism $\Omega/N\ra\Omega/\Gamma$ is finite topological, hence it is finite \'etale.  Because Mumford curves are analytifications of
proper schemes, and the analytification functor for proper schemes is fully faithful by theorem \ref{curveana}, the associated morphism of
algebraic curves is finite \'etale of degree $|\Gamma/N|$; this is proved in \cite[Proposition 3.3.11]{Berk1993}.
\end{proof}

\subsubsection{Examples}\label{subsexamples}\text{ }\\
We give examples of elliptic and hyperelliptic Mumford curves.
\begin{defin}
Let $q\in\K$ with $0<\n{q}<1$, set $\Sigma:=\{0,\infty\}$ and let $\Gamma$ be the cyclic group that is generated by $\begin{pmatrix} q & 0 \\ 0 & 1 \end{pmatrix}$. The group $\Gamma$ is a Schottky group and acts freely on $\PP(\Cp)-\Sigma=\Gm(\Cp)$. The associated $\K$-Mumford curve is called {\em Tate curve}.
\end{defin}
\begin{nota}
Let $\K$ be a local number field and assume that its residue field does not have characteristic 2. 
In \cite[page 167]{Gerr1980} the possible stable reductions of a $\K$-curve of genus 1 are calculated. It follows that for an algebraic $\K$-curve $X$ the following properties are equivalent:
\begin{enumerate}
\item It is a Tate curve.
\item It is an elliptic curve that has a stable $\K^{\circ}$-model with bad reduction.
\end{enumerate}
\end{nota}

Let $\K$ be a non-Archimedean field.
\begin{defin}
Given a family of elements $(s_{i})_{i=0,...,g}$ of order 2 in the group $\PGL_{2}(\K)$ such that the group $\Gamma:=\langle(s_{i})_{i=0,...,g}\rangle$ is discontinuous and isomorphic to the free product of the groups $\{1,s_{i}\}$ for $i=0,...,g$. Let $W$ be the kernel of the group morphism defined by $s_{i}\mapsto -1$ for $i=0,...,g$. The group $W$ is called {\em Whittaker group}.
\end{defin}
\begin{nota}
The group $W$ is a free group generated by $s_{1}s_{0},...,s_{g}s_{0}$ and it is discontinuous. Thus it is a Schottky group.
\end{nota}
\begin{defin}
A proper and smooth algebraic curve $X/\K$ is called {\em hyperelliptic curve} if there is a morphism $X\ra\PP$ of degree 2.
\end{defin}

\begin{proposition}[\cite{Gerr1980}, page 282]
For an algebraic curve $X$ over $\K$ the following properties are equivalent:
\begin{enumerate}
\item $X$ is a Mumford curve, and there is a Whittaker group $W$ such that $X^{an}=\Omega_{W}/W$.
\item The curve $X$ is hyperelliptic and has a stable $\widetilde{\K}$-split degenerated model.
\end{enumerate}
\end{proposition}

\subsubsection{The Reduction graph of a Mumford curve}\label{subsred}\text{ }\\
In this paragraph the ground field $\K$ is assumed to have discrete valuation.\\
In the articles \cite{Falt1983} and \cite{Put1986} by Faltings and van der Put--Reversat it is important to consider the graph associated to a specific reduction of $X^{rig}$. This reduction does not correspond to the stable model $\X^{stab}$ of $X$ but to the minimal resolution $\X^{min}$ of $\X^{stab}$. Refering to page 163 of Mumford's article \cite{Mumf1972} we can define a certain reduction on $\Omega$ such that it becomes a formal rigid space and the quotient $\Omega^{form}/\Gamma$ of the associated formal scheme by the Schottky group $\Gamma$ is isomorphic to the formal completion of $\X^{min}$.\\
The minimal resolution of $\X^{stab}$ exists, because $\X^{stab}$ is an excellent, reduced and Noetherian scheme of dimension 2 as a result of Lipman's theorem of resolution of singularities \cite[Theorems 8.2.39 and 8.3.44]{Liu2002}. The scheme $\X^{min}$ is projective by the result \cite[Theorem 8.3.16]{Liu2002} of Lichtenbaum. It is the minimal regular model, because of theorem \cite[Theorem 9.3.21]{Liu2002} and the fact that the arithmetic genus of the generic fibre is greater than 0, it is unique up to unique isomorphism. The scheme $\X^{min}$ is not stable anymore but only semistable. Its special fibre has only rational components, all its double points are $\widetilde{\K}$-rational and every non-singular component meets the other components in at least 2 points. From now on the reductions attached to $\Omega$ and $X^{rig}$ are those whose associated model $\X$ is the minimal regular model of $X$.
\begin{defin}

Let $\widetilde{X}$ be a $\widetilde{\K}$-reduction of a rigid space. Assume that $\widetilde{X}$ is reduced, connected, one-dimensional and locally of finite type, having at most ordinary double points and these to be $\widetilde{\K}$-rational, and assume that all its components are rational. Define a graph by one-to-one correspondence between irreducible components and {\em vertices}, and between double points and {\em edges} such that a component contains a double point if and only if the corresponding vertex is an endpoint of the corresponding edge. The graph obtained is called {\em graph of the reduction}. We express it by the term $T(\widetilde{X})$.   
\end{defin}
\begin{nota}
The reduction graph $T(\widetilde{\Omega})$ of the universal covering of a Mumford curve $X^{rig}:=\Omega/\Gamma$ is a tree, on which $\Gamma$ operates freely. It is true that \[T(\widetilde{X^{rig}})=T(\widetilde{\Omega})/\Gamma\] and the graph $T(\widetilde{X^{rig}})$ is finite.
\end{nota}
\section[Various Constructions]{Representations attached to vector bundles}
\subsection{A review of the constructions}
In this subsection we present various constructions of representations of fundamental groups. The construction by van der Put and Reversat in \cite{Put1986} is for Mumford curves over non-Archimedean fields and uses rigid geometry. For every semistable vector bundle of degree zero a representation of the topological fundamental group is defined. This work is a generalisation of the article \cite{Falt1983} by Faltings. The most recent contribution in this direction are the articles \cite{Deni2003}, \cite{Deni2004} and \cite{Deni2004b} by Deninger and Werner. Given a projective and smooth curve, for a certain class of vector bundles on it they construct continuous representations of the curve's algebraic fundamental group.

\subsubsection{The construction of van der Put and Reversat}\text{ }\\
This paragraph refers to the article \cite{Put1986}. We describe the construction only in the restricted situation that the considered Mumford curve comes from a local number field. Let $X'^{rig}:=\Omega'/\Gamma$ be a Mumford curve over a local number field $\K'$. Denote the analytic reductions of $\Omega'$ and $X'^{rig}$ by $R':\Omega'\ra\widetilde{\Omega'}$ and $r':X'^{rig}\ra\widetilde{X'^{rig}}$ respectively. We remind the reader that the reduction $r':X'^{rig}\ra\widetilde{X'^{rig}}$ is associated to the minimal regular model of $X'$ as explained in section \ref{chpre} in the paragraph \ref{subsred} on Mumford curves. Let $\K$ be a complete subfield of $\Cp$, and let $X, \Omega,  r$ and $R$ be the respective base extension to $\K$. We obtain the following commutative diagram.
\[\xymatrix{\Omega \ar[r]^{u}\ar[d]^{R} & X^{rig} \ar[d]^{r} 
\\{\widetilde{\Omega} \ar[r]^{\widetilde{u}}} & {\widetilde{X^{rig}}}}\]
Given an orientation on the reduction graph of $X^{rig}$, one obtains a $\Gamma$-invariant orientation on the reduction graph of $\Omega$. 

If $d$ is a double point of $\widetilde{\Omega}$, then there is a finite extension field $\L$ of $\K$ such that
\[R^{-1}(d)\cong\{x\in\Aff\suchthat \n{\pi_{d}}<T(x)<1\}\text{ and }\n{\pi_{d}}\in \L^{*}.\] We assume that $\K$ equals $\L$. The general theorem follows by finite Galois descent.
The number $\n{\pi_{d}}$ is invariant under the action of $\Gamma$, i. e. $\n{\pi_{d}}=\n{\pi_{\gamma(d)}}$.

\begin{defin}
Let $L_{1}, L_{2}$ be irreducible components of $\widetilde{\Omega}$ and define 
\[d_{+}(L_{1},L_{2}):=\prod_{+}\n{\pi_{d}},\]
where $\Pi_{+}$ denotes the product over those edges on a geodesic between $L_{1}$ and $L_{2}$ whose orientation is the same as the orientation of the path from $L_{1}$ to $L_{2}$.
\end{defin}
\begin{nota}
Let $L$ be an irreducible component of $\widetilde{\Omega}$. Since the number of edges on a geodesic between $L$ and $\gamma L$ is at least two, there is an orientation of the reduction graph of $X^{rig}$ such that the induced $\Gamma$-invariant orientation of the reduction graph of $\Omega$ has the property 
\[d_{+}(L,\gamma L)\ne 1.\]
This can be achieved by orienting at least one edge between $L$ and $\gamma L$ positively and at least one edge negatively.\\
In the following we will always assume such an orientation.
\end{nota}
\begin{defin}
Let $L_{0}$ be an irreducible component of $\widetilde{\Omega}$ and $V$ a $\K$-vector space of finite dimension with fixed basis $B$. Given a matrix $\sigma$ in $\Gl(V)$ with entries $\sigma_{ij}$ with respect to the basis $B$, we define \[\n{\sigma}:=\max_{i,j=1,\ldots,n} \n{\sigma_{ij}}.\]
A representation $\rho:\Gamma\ra\Gl(V)$ is called {\em $\phi$-bounded} if
\begin{enumerate}
\item $\sup_{\gamma\in\Gamma}(d_{+}(L_{0},\gamma L_{0})\n{\rho(\gamma^{-1})})$ is finite and
\item $\{\gamma\in\Gamma \suchthat d_{+}(L_{0},\gamma L_{0})\n{\rho(\gamma^{-1})}\geq \epsilon\}$ is well-ordered for every $\epsilon > 0$.
\end{enumerate}
\end{defin}
\begin{nota}
The definition does not depend on the choice of $B$ or $L_{0}$.
\end{nota}
\begin{lemma}
Every representation $\rho$ with the property $\n{\rho(\gamma)}=1$ for all $\gamma\in\Gamma$ is $\phi$-bounded.
\end{lemma}
\begin{proof}
The first condition of $\phi$-boundedness is satisfied because \[d_{+}(L_{0},\gamma L_{0})<1.\] To prove the second property
we first note that there are only finitely many different values of $\n{\pi_{d}}$ since this number is $\Gamma$-invariant and $T(\widetilde{\Omega})/\Gamma$ is finite. Let $q<1$ be the biggest of these numbers. If $\length_{+}(L_{0},\gamma L_{0})$ denotes the number of positively oriented edges and $\length(L_{0},\gamma L_{0})$ the number of all edges on a geodesic between $L_{0}$ and $\gamma L_{0}$, then there is a number $C$ such that 
\[d_{+}(L_{0},\gamma L_{0})\geq \epsilon \Ra q^{l_{+}(L_{0},\gamma L_{0})}\geq \epsilon \Ra l_{+}(L_{0},\gamma L_{0})\leq C.\]  
Because of the chosen orientation, $l_{+}(L_{0},\gamma L_{0})\geq \frac{1}{n}l(L_{0},\gamma L_{0})$, where $n$ is the length of a geodesic of maximal length in the graph $T(\widetilde{X^{rig}})$. Hence it follows
\[ d_{+}(L_{0},\gamma L_{0})\geq \epsilon \Ra \frac{1}{n}l(L_{0},\gamma L_{0})\leq C \Ra l(L_{0},\gamma L_{0}) \leq C\cdot n.\]
Only a finite number of group elements $\gamma$ in $\Gamma$ can satisfy this, since the reduction graph of $\Omega$ is a locally finite tree.\\
Thus the set \[\{\gamma\in\Gamma \suchthat d_{+}(L_{0},\gamma L_{0})\n{\rho(\gamma^{-1})}\geq \epsilon\}\]
is finite and therefore well-ordered.
\end{proof}

\begin{defin}\label{rvdpdefin}\text{ }
\begin{enumerate}
\item Given a $\K$-vector space $V$ with basis $e_{1},...,e_{r}$ and a representation $\rho:\Gamma\ra\Gl(V)$, define the free $\Loc_{\Omega}(\Omega)$-module 
\[W:=\Loc_{\Omega}(\Omega)\otimes_{\K}V\]
endowed with the $\Gamma$-action 
\[\gamma(\sum_{i=1}^{r} f_{i}\otimes e_{i}):=\sum_{i=1}^{r}(f_{i}\circ\gamma^{-1})\otimes\rho(\gamma)e_{i}.\]
\item  Define a vector bundle $E_{\rho}$ by \[E_{\rho}(U):=\{f\in W\otimes_{\Loc_{\Omega}(\Omega)}\Loc_{\Omega}(u^{-1}(U))\suchthat\gamma(f)=f, \text{for all }\gamma\in\Gamma\}\] for every open affinoid domain $U\subset X$;  here the $\Gamma$-action is defined by 
\[\gamma(\sum_{i=1}^{r} e_{i}\otimes f_{i}):=\sum_{i=1}^{r}\gamma(e_{i})\otimes (f_{i}\circ\gamma^{-1}).\] Below we use the shorter notation 
\[E_{\rho}=[V\otimes u_{*}\Loc_{\Omega}]^{\Gamma}.\]
\end{enumerate}
\end{defin}
\newpage
\begin{theorem}[\cite{Put1986}, Th\'eor\`eme principal 5]\text{ }
\begin{enumerate}
\item Every semistable vector bundle $E$ of rank $\rk$ and degree 0 on a Mumford curve is isomorphic to $E_{\rho}$ for an $\rk$-dimensional $\phi$-bounded representation $\rho:\Gamma\ra\Gl(V)$.
\item For any $\phi$-bounded $\rk$-dimensional representation $\rho:\Gamma\ra\Gl(V)$ the vector bundle $E_{\rho}$ is semistable of degree 0 and has rank $\rk$.
\item If $Hom_{\Gamma}(\rho_{1},\rho_{2})$ denotes homomorphisms of $\Gamma$-representations, then there is a natural isomorphism \[\Hom_{\Loc_{X}}(E_{\rho_{1}},E_{\rho_{2}})\cong Hom_{\Gamma}(\rho_{1},\rho_{2}).\]
\end{enumerate}
\end{theorem}
\begin{corollar}
The functor $\rho\mapsto E(\rho)$ is an equivalence between the category of $\phi$-bounded representations and the category of semistable vector bundles of degree zero.   
\end{corollar}
\begin{proof}
We have to prove that the functor is fully faithful and essentially surjective. This follows immediately from the above theorem.
\end{proof}
\begin{nota}\label{prremarks}\text{ }
\begin{enumerate}
\item As an abbreviation we call a van der Put--Reversat representation a {\em PR-representation}. 
\item If $E$ is a semistable vector bundle of degree 0, then the attached representation is described explicitly by \cite[5.8, 5.9]{Put1986} as the natural $\Gamma$-action on a sub $\K$-vector space of $(\widetilde{u}^{*}r_{*}E)(\widetilde{\Omega})=u^{*}E(\Omega)$. Furthermore with the notations in definition \ref{rvdpdefin} it is true that $W\cong u^{*}E(\Omega)$.
\item The above theorem gives an equivalence of categories between semistable vector bundles of degree 0 and $\phi$-bounded representations. Unfortunately this equivalence does not commute with tensor products and duals (see remark of Faltings in \cite{Falt1983}). This is clear even in the line bundle case, since taking the tensor product of two line bundles is associated to the multiplication in the Jacobian and taking the dual is associated to taking the inverse in the Jacobian, but the $\phi$-bounded representations are not closed under this operations.
\end{enumerate}
\end{nota}
We are going to characterise those vector bundles whose representation has image in $\Gl_{\rk}(\K^{\circ})$. First we define two categories of vector bundles.
\begin{defin}
Let $\K$ be a complete subfield of $\Cp$. Given a $\K$-Mumford curve that comes from a Mumford curve $X'$ over a local number field $\K'$, let $i':X'\inj \X'$ be its uniquely defined minimal regular $\K'^{\circ}$-model. Let $\X$ be $\X'\otimes_{\K'^{\circ}}\K^{\circ}$ and denote the canonical morphism by $i:\X\inj X$.
\begin{enumerate}
\item Define $\BB_{X}^{\X}$ as the full subcategory of all vector bundles on $X$ whose objects have the following properties: they are semistable of degree zero, and for every $E\in\BB_{X}^{\X}$ there is a vector bundle $\E$ on $\X$ with $i^{*}\E\cong E$.
\item Define $\BB_{X^{rig}}^{\widehat{\X}}$ as the full subcategory of all semistable rigid vector bundles of degree zero on $X^{rig}$ whose associated van der Put--Reversat representation is isomorphic to a representation that has image in $\Gl_{\rk}(\K^{\circ})$.
\end{enumerate}
\end{defin}
\begin{proposition}\label{modelandrep}
Rigidication of vector bundles induces an equivalence of categories between $\BB_{X}^{\X}$ and $\BB_{X^{rig}}^{\widehat{\X}}$.
\end{proposition}
\begin{proof}
We have to show that the rigidification functor maps elements of $\BB_{X}^{\X}$ into $\BB_{X^{rig}}^{\widehat{\X}}$ and that algebraisation of rigid vector bundles maps elements of $\BB_{X^{rig}}^{\widehat{\X}}$ into $\BB_{X}^{\X}$.
\begin{enumerate}
\item $\BB_{X}^{\X}\ra\BB_{X^{rig}}^{\widehat{\X}}$\\
If $E'$ is a vector bundle in $\BB_{X}^{\X}$, then it is isomorphic to $E:=i^{*}\E$ for a vector bundle $\E$ on $\X$. Let $\rho:\Gamma\ra\Gl(V)$ be the PR-representation that is associated to $E^{rig}$. As remarked in \ref{prremarks} the space $V$ is a sub $\K$-vector space of $(\widetilde{u}^{*}r_{*}E^{rig})(\widetilde{\Omega})$. The action by $\Gamma$ is $\gamma f:=\rho(\gamma)(f\circ\gamma^{-1})$. Denote by $u^{form}$ the morphism \[u^{form}:(\widetilde{\Omega},R_{*}\Loc_{\Omega}^{\circ})\ra(\widetilde{X^{rig}},r_{*}\Loc_{X^{rig}}^{\circ})\] between formal schemes. Notice that on the underlying topological spaces the morphisms $u^{form}$ and $\widetilde{u}$ coincide.\\ 
If $\widehat{\E}$ denotes the formal completion of $\E$, then it follows 
\[
\begin{array}{rrl}
\widetilde{u}^{*}r_{*}E^{rig} & \cong &\widetilde{u}^{*}(\widehat{\E}\otimes r_{*}\Loc_{X^{rig}})\\
\comment{&\cong&\widetilde{u}^{-1}(\widehat{\E}\otimes r_{*}\Loc_{X^{rig}})\otimes R_{*}\Loc_{\Omega}\\
&\cong&\widetilde{u}^{-1}\widehat{\E}\otimes R_{*}\Loc^{\circ}_{\Omega}\otimes \widetilde{u}^{-1}r_{*}\Loc_{X^{rig}}\otimes R_{*}\Loc_{\Omega}\\}
&\cong&(u^{form})^{*}\widehat{\E}\otimes \widetilde{u}^{*}r_{*}\Loc_{X^{rig}}.
\end{array}
\]
Therefore up to natural isomorpism the $\K$-vector space $V$ is a subspace of 
\[((u^{form})^{*}\widehat{\E}\otimes \widetilde{u}^{*}r_{*}\Loc_{X^{rig}})(\widetilde{\Omega}).\]
The action of $\Gamma$ on $\widetilde{u}^{*}r_{*}\Loc_{X^{rig}}(\widetilde{\Omega})$ is trivial. Therefore $\Gamma$ acts only on $(u^{form})^{*}\widehat{\E}$, hence $\rho(\gamma)\in\Gl_{\rk}(\K^{\circ})$ for all $\gamma\in\Gamma$.
\item $\BB_{X^{rig}}^{\widehat{\X}}\ra\BB_{X}^{\X}$\\
Each vector bundle in $\BB_{X^{rig}}^{\widehat{\X}}$ is the rigidification of an algebraic vector bundle $E$ that is semistable and of degree zero. Since $\Bild\rho\subset\Gl_{\rk}(\Ok)$, the bundle $E^{rig}$ has the formal model $\E$ defined by \[\E:=[V^{\circ}_{\rho}\otimes_{\K^{\circ}} r_{*}u_{*}\Ooo]^{\Gamma},\] where $V^{\circ}$ is defined as follows. If $V=\oplus_{i=1}^{\rk}e_{i}\K$, then define 
\[V^{\circ}:=\oplus_{i=1}^{\rk}e_{i}\K^{\circ}.\]
By corollary \ref{formelvectgaga} the bundle $\E$ is the formal completion of an algebraic vector bundle $\E'$ on $\X$. By lemma \ref{modelgaga} the vector bundle $\E'$ is an algebraic $\X$-model of $E$.
\end{enumerate}
\end{proof}
\begin{defin}\label{canmodel}
Let $E$ be a vector bundle in $\BB_{X}^{\X}$. By the proposition above there is a vector bundle model $\E$ of $E$ on $\X$ with the property \[\widehat{\E}=[V^{\circ}_{\rho}\otimes_{\K^{\circ}} r_{*}u_{*}\Ooo]^{\Gamma}.\] The bundle $\E$ is called the {\em canonical model} of $E$. 
\end{defin}
\subsubsection{The construction of Faltings}\text{ }\\
Let $X$ be a Mumford curve over a discrete non-Archimedean field $\K$. In his article \cite{Falt1983} Faltings introduced the concept of $\phi$-boundedness and proved for the first time the correspondence between $\phi$-bounded representations of $\pi^{top}_{1}(X^{an},\overline{x})$ and semistable vector bundles of degree zero. Because in his proof he used the theory of formal schemes, he was limited to the case of discrete valuation. In what follows we will solely deal with the van der Put--Reversat construction, therefore we are not going to go into the details of Faltings' article, but we only remark that in the case of a discrete field both constructions attach isomorphic $r_{*}u_{*}\Loc_{\Omega}$-vector bundles to a $\phi$-bounded representation in the following sense: since $\phi$-boundedness depends on some choices (e.g. the orientation of the reduction graph) we assume that in Faltings' article and in the article of van der Put and Reversat the same choices are made, such that the same representations are $\phi$-bounded.\\ 
Let $u^{form}:\Omega^{form}\ra X^{form}$ be the quotient morphism between the formal schemes that are associated to the minimal regular model of $X$ (for details see section 1 paragraph \ref{subsred}). In Faltings' article to a $\phi$-bounded representation $\rho$ the vector bundle
\[[(\K\otimes_{\K^{\circ}}u^{form}_{*}\Loc_{\Omega^{form}})^{\rk}]^{\Gamma}\] is attached. Because of $\Loc_{\Omega^{form}}\otimes\K\cong r_{*}\Loc_{\Omega}$ this gives the vector bundle 
\[[\K^{\rk}\otimes_{\K}(\widetilde{u}_{*}R_{*}\Loc_{\Omega})]^{\Gamma}.\] It is naturally isomorphic to \[r_{*}[\K^{\rk}\otimes_{\K}(u_{*}\Loc_{\Omega})]^{\Gamma}.\] 
This vector bundle defines uniquely the rigid vector bundle \[E:=[\K^{\rk}\otimes_{\K}(u_{*}\Loc_{\Omega})]^{\Gamma}.\]
$E$ is attached to $\rho$ by the construction of van der Put and Reversat.
\subsubsection{The construction of Deninger and Werner}\text{ }\\
In this paragraph we review the constructions that are made in \cite{Deni2004}. Let $R$ be a valuation ring with quotient field $Q$ of characteristic 0, and let $X$ be a projective and smooth curve over $Q$ with an $R$-model $\X$.
To begin with, we define certain categories of coverings of $\X$.
\begin{defin}\text{ }
\begin{enumerate}
\item The objects of the category $\Cov^{fpp}(\X)$ are the finitely presented proper $R$-morphisms $\pi:\Y\ra\X$ that satisfy 
$\pi\otimes Q \in \Cov^{alg}(X)$. The morphisms are commutative triangles
\[
\xymatrix{{\Y_{1}} \ar[rr] \ar[rd] && {\Y_{2}} \ar[ld]\\
&{\X}.}\]
\item The category $\Cov^{\text{good}}(\X)$ is the full subcategory of $\Cov^{fpp}(\X)$ consisting of those morphisms $\pi:\Y\ra\X$ whose structure morphism $\lambda:\Y\ra\spec R$ is flat, the equality $\lambda_{*}\Loc_{\Y}=\Loc_{\spec R}$ holds universally and the morphism $\Y_{Q}\ra\spec Q$ is smooth.
\end{enumerate}
\end{defin}
Let $X$ be a projective and smooth curve over $\overline{\Qp}$ and $\X$ a $\overline{\Zp}$-model of $X$.
Let $\VVec_{\X_{\op}}$ be the category of all vector bundles on $\X_{\op}$. Denote by $\cdot_{n}$ the reduction modulo $p^{n}$ and define the following two categories of vector bundles.
\begin{defin}\text{ }
\begin{enumerate}
\item Define $\BB_{\X_{\op}}$ as the full subcategory of $\VVec_{\X_{\op}}$ whose objects $\E$ have the property that for every natural number $n$ there is a covering $\pi:\Y\ra\X$ in $\Cov^{fpp}(\X)$, such that
\[\pi^{*}_{n}\E_{n}\cong\Loc_{\Y_{n}}^{\rk \E}.\] 
\item The category $\BB_{X_{\Cp}}$ is the full subcategory of $\VVec_{X_{\Cp}}$ consisting of all vector bundles $E$ on  $X_{\Cp}$ 
that have the following property: it exists a $\overline{\Zp}$-model $i:X\inj \X$ of $X$ and a vector bundle $\E$ in $\BB_{\X_{\op}}$, such that $E$ is isomorphic to $i^{*}\E$.
\end{enumerate}
\end{defin}

The construction of Deninger and Werner, as given in \cite{Deni2004}, is as follows: Let $\E$ be a bundle in $\BB_{\X_{\op}}$ and $\overline{x}$ a geometric 
point in $X(\Cp)$; by properness this gives $x_{\op}\in\X(\op)$ and by reduction modulo $p^{n}$ we obtain a section
\[\xymatrix{{x_{n}: \spec \on} \ar[r]& {\spec \op} \ar[r]^{x_{\op}} & {\X_{\op}}.}\]
We write $\E_{x_{\op}}$ for $x_{\op}^{*}\E$ viewed as a free $\op$-module of rank $\Rang \E$, and define \[\E_{x_{n}}:=x_{\op}^{*}\E\otimes\on\] viewed 
as a free $\on$-module of rank $\Rang \E$. The module $\E_{x_{n}}$ is a topological $\on$-module endowed with the discret topology. It is true that \[\plim_{n\ra \infty}\E_{x_{n}}=\E_{x_{\op}}\] as topological modules.

Denote by $\FMod_{\op}$ the category of free topological $\op$-modules of finite rank and define a continuous functor
\[\rho_{\E}:\Pi_{1}(X_{\Cp})\ra \FMod_{\op}.\]
On objects it is defined by 
\[\rho_{\E}(\overline{x}):=\E_{x_{\op}},\]
and on homomorphisms the map \[\rho_{\E}:\Iso(F_{\overline{x}},F_{\overline{x'}})\ra\Hom_{\op}(\E_{x_{\op}},\E_{x_{\op}'})\] is defined as the projective limit of 
\[\rho_{\E,n}:\Iso(F_{\overline{x}},F_{\overline{x'}})\ra\Hom_{\on}(\E_{x_{n}},\E_{x_{n}'}).\]\\
For every natural number $n$ the functor $\rho_{\E,n}$ is defined as follows. By \cite[Corollary 3 3)]{Deni2004} it is true that there is a covering $\pi:\Y\ra\X$ in $\Cov^{\text{good}}(\X)$, such that
\[\pi^{*}_{n}\E_{n}\cong\Loc_{\Y_{n}}^{\rk\E}.\]  
Set $Y:=\Y\otimes\overline{\Q_{p}}$, then $Y\ra X$ is a finite \'etale covering. Choose $\overline{y}$ above $\overline{x}$. A path $\gamma$ from $\overline{x}$ to $\overline{x'}$ is an isomorphism of fibre functors $F_{\overline{x}}$ to $F_{\overline{x'}}$. Define $\overline{y'}:=\gamma \overline{y}$; the point $\overline{y'}$ lies over $\overline{x'}$. Because 
$\lambda_{*}\Loc_{\Y}=\Loc_{\spec \overline{\Zp}}$ holds universally, it follows $\lambda_{n*}\Loc_{\Y_{n}}=\Loc_{\spec \on}.$ Therefore
\[
\xymatrix{y^{*}_{n}: \Gamma(\Y_{n},\pi^{*}_{n}\E_{n}) \ar[r] & {\Gamma(\spec \on, y_{n}^{*}\pi_{n}^{*}\E_{n}) =  \E_{x_{n}}}}
\]
is an isomorphism. This follows, because $\Gamma(\Y_{n},\Loc_{\Y_{n}})\cong \on$. Define the morphism \[\rho_{\E,n}(\gamma):=(\gamma y_{n})^{*}\circ (y_{n}^{*})^{-1}: \E_{x_{n}}\ra \E_{x'_{n}}.\]
Define 
\[\rho_{\E}:=\plim_{n}\rho_{\E,n}.\]
By construction $\rho_{\E}$ is continuous. It is independent of all choices.\\
The last step is the definition of the functor \[\rho:\BB_{\X_{\op}}\ra\Rep_{\Pi_{1}(X)}(\op).\] On objects it is defined by \[\rho(\E):=\rho_{\E},\] and if $f:\E_{1}\ra \E_{2}$ is a morphism of vector bundles, then define
\[\rho(f):=\{\E_{1,x_{\op}}\ra \E_{2,x_{\op}} | x\in\Ob\Pi_{1}(X)\}.\]
This set is a natural transformation from $\rho_{\E_{1}}$ to $\rho_{\E_{2}}$.

\begin{nota}
If $X$ is defined over $\K$, then the group $\Gal(\overline{\Qp}/\K)$ acts from the left on the categories $\BB_{\X_{\op}}$ and $\Rep_{\Pi_{1}(X)}(\op)$ resp. on the categories $\BB_{X_{\Cp}}$ and $\Rep_{\Pi_{1}(X)}(\Cp)$. The action is defined on page 39 of \cite{Deni2004}. Because we do not need the action, we do not give its definition.
\end{nota}
\begin{theorem}[\cite{Deni2004}, Proposition 24]
The functor \[\rho:\BB_{\X_{\op}}\ra\Rep_{\Pi_{1}(X)}(\op)\text{: }\rho(\E):=\rho_{\E}\] is $\op$-linear and exact, commutes with duals, tensors, internal homs and exterior powers of vector bundles. Exact sequences are mapped to exact sequences of representations. If $X$ is defined over $\K$, then the functor commutes with the left action of $\Gal(\overline{\Qp}/\K)$ on the categories $\BB_{\X_{\op}}$ and $\Rep_{\Pi_{1}(X)}(\op)$. 
\end{theorem}

The construction can be carried over to $\BB_{X_{\Cp}}$. If $E$ is an object of $\BB_{X_{\Cp}}$, then by definition $E\cong j_{\X_{\op}}^{*}\E$. Therefore
\[\xymatrix{
{\psi_{\overline{x}}= \overline{x}^{*}\psi: E_{\overline{x}}} \ar[r]^-{\cong} & (j_{\X_{\op}}^{*}\E)_{\overline{x}}=\E_{x_{\op}}\otimes\Cp}.
\]
Define $\rho_{E}:\Pi_{1}(X)\ra\FMod_{\Cp}$ as 
\[\rho_{E}(\overline{x}):=\overline{x}^{*}E=E_{\overline{x}}\text{ and  }\rho_{E}(\gamma)=\psi^{-1}_{\overline{x}}(\rho_{\E}(\gamma)\otimes\Cp)\psi_{\overline{x}}.\]

\begin{theorem}[\cite{Deni2004}, Theorem 28]\label{dewecp}
The above defined functor $\rho$ is $\Cp$-linear and exact, commutes with duals, tensors, internal homs and exterior powers of vector bundles. If $X$ is defined over $\K$, then the functor commutes with the left action of $\Gal(\overline{\Qp}/\K)$ on the categories $\BB_{X_{\Cp}}$ and $\Rep_{\Pi_{1}(X)}(\Cp)$.
\end{theorem}
\begin{nota}
We call the Deninger--Werner representation also {\em DW-re\-pre\-sen\-tation}.
\end{nota}

The Deninger--Werner construction can be performed for a more general category $\C_{X_{\Cp}}$ of vector bundles, which contains $\BB_{X_{\Cp}}$. These vector bundles can be characterised by the following definition.
\begin{defin}
\begin{enumerate}
\item A vector bundle $E$ on a smooth projective curve over a field of characteristic $p>0$ is called {\em strongly semistable} if the pullbacks of $E$ by all non-negative powers of the absolute Frobenius are semistable.
\item Let $R$ be a valuation ring with quotient field $Q$ and residue field $k$. Let $X$ be a smooth projective curve over $Q$ with $R$-model $\X$. A vector bundle $\E$ on $\X$ is called {\em strongly semistable} if the pullback of $\E_{k}$ to the normalisation of each irreducible component of $\X_{k}$ is strongly semistable.
\item Let $X$ be a smooth projective curve over $\overline{\Qp}$. A vector bundle $E$ on $X_{\Cp}$ is said to have {\em strongly semistable reduction} if there is a $\overline{\Z_{p}}$-model $\X$ of $X$ such that $E$ extends to a vector bundle on $\X$ that is strongly semistable.
\item A bundle $E$ on $X$ is said to have {\em potentially strongly semistable reduction} if there is a finite \'etale morphism $\pi:Y\ra X$ of smooth projective curves sucht that $\pi^{*}E$ has strongly semistable reduction.
\end{enumerate}
\end{defin}
The authors prove that $\C_{X_{\Cp}}$ is exactly the category of vector bundles of degree zero with potentially strongly semistable reduction. The \'etale parallel transport and in particular the associated respresentation of the algebraic fundamental group exist in this case. The parallel transport is compatible with tensor products, duals, internal homs, pullbacks and Galois conjugation (cf. \cite[Theorem, page 2]{Deni2004}). 
\subsection{Comparison of the constructions}
In this subsection the construction of van der Put and Reversat will be compared to the one of Deninger and Werner. This will be done separately for vector bundles that are defined over a local number field (the case of discrete valuation) and for vector bundles that are only defined over $\Cp$ (the general case).
\subsubsection{Descent of vector bundle reductions}\text{ }\\
This paragraph serves as a preparation for the general case and shows that it is possible to perform the construction of the DW-representation modulo $p^{n}$ over a discrete valuation ring instead of over the ring of integers $\op$ of $\Cp$\\
Let $S_{1}$ be a scheme and $S_{0}$ a closed subscheme of $S_{1}$ that is defined by an ideal $J$ of square zero. Let $G$ be an $S_{1}$-group scheme. For every $S_{1}$-scheme $X$ write $X_{0}$ for the base change $X\times_{S_{1}}S_{0}$. If $P_{0}$ is a $G_{0}$-torsor over $S_{0}$ for the Zariski (resp. \'etale) topology and $\Lie(G_{0}/S_{0})'$ is the $\Loc_{S_{0}}$-module obtained by twisting $\Lie(G_{0}/S_{0})$ by $P_{0}$, then define \[W:=\Lie(G_{0}/S_{0})'\otimes_{\Loc_{S_{0}}}J.\] 
With this notations the following holds:
\begin{theorem}[\cite{Gira1971}, VII Th\'eor\`eme 1.3.1]
If $G$ is smooth, then it exists a class $c\in H^{2}_{\text{Zar}}(S_{0},W)$  whose vanishing is necessary and sufficient for the existence of a $G$-torsor $P_{1}$ for the Zariski (resp. \'etale) topology on $S_{1}$ that lifts $P_{0}$. If $c=0$, then the set of isomorphism classes of liftings of $P_{0}$ is a $  H^{1}_{\text{Zar}}(S_{0},W)$-torsor.
\end{theorem}
\begin{lemma}\label{reduction}
Let $\K$ be a local number field with ring of integers $\K^{\circ}$ and $\X$ a projective flat $\K^{\circ}$-curve. Let $\op$ be the ring of integers of $\Cp$ and let $\E$ be a vector bundle on $\X_{\op}:=\X\otimes\op$, then for every natural number $n$ there is a local number field $\L^{(n)}/\K$ and a vector bundle $\E^{(n)}$ on $\X_{\Loc_{\L^{(n)}}}$, such that \[\E\otimes \op/p^{n}=\E^{(n)}\otimes \Loc_{\L^{(n)}}/p^{n}.\]
\end{lemma}
\begin{proof}
We use the notations in the theorem. In our case we consider only the Zariski topology and the smooth $S_{1}$-group scheme $G:=\Gl_{\rk}$.
Then $P_{0}$ is a vector bundle on $S_{0}$ with respect to the Zariski topology. Because of $\op/p^{n}=\overline{\Zp}/p^{n}$ it follows \[H^{1}(\X\otimes\op/p^{n},\Gl_{\rk})=H^{1}(\X\otimes\overline{\Zp}/p^{n},\Gl_{\rk}),\]
and the vector bundle $P_{0}$ is already defined on $\X_{\Loc_{\L^{(n)}}}$ for a finite extension field $\L^{(n)}$ of $\K$. 
We use the theorem in the following situation. Define 
\[\begin{array}{l}
S_{i}:=\X_{n+i}:=\X\otimes\Loc_{\L^{(n)}}/p^{n+i}\\
 J:=p^{n+i-1}(\Loc_{\L^{(n)}}/p^{n+i}) \text{ and } G:=\Gl_{\rk}.
\end{array}\]
We proceed by induction. The induction start is $i=1$.\\ 
Since $\Lie(G_{0}/S_{0})'$ is a locally free $\Loc_{S_{0}}$-module (cf. \cite[VII 1.3.1.1]{Gira1971}), the sheaf $W$ is a coherent sheaf on $S_{0}$. The cohomology group $H^{2}_{\text{Zar}}(S_{0},W)$ vanishes, since $S_{0}$ is a curve and $W$ is coherent. Therefore by the theorem it exists a vector bundle $\E_{1}^{(n)}:=P_{1}$ on $S_{1}$ with $\E_{1}^{(n)}\otimes\Loc_{\L^{(n)}}/p^{n}=P_{0}$.\\
We have a chain of closed subschemes 
\[S_{0}\inj S_{1}\inj \ldots \inj S_{i} \inj S_{i+1} \inj \ldots.\]
In the induction step we apply the theorem to $S_{i}$ and $S_{i+1}$ and we get a projective system $(\X_{n+i},\E^{(n)}_{i})_{i\geq 1}$. The projective limit \[\plim_{i\geq 1}\E^{(n)}_{i}\] is a formal vector bundle on the formal scheme 
\[\ilim_{i\geq n+1}\X_{i}=\widehat{\X\otimes\Loc_{\L}}.\] By Grothendiecks GAGA \ref{formelvectgaga} there is an algebraisation $\E^{(n)}$ of $\plim\E^{(n)}_{i}$, which is a vector bundle. It follows
$\E^{(n)}\otimes\Loc_{\L^{(n)}}/p^{n}=\E\otimes\op/p^{n}$.   
\end{proof}

\subsubsection{The case of discrete valuation}\text{ }\\
Let $\K$ be a local number field and $X$ a Mumford curve over $\K$. Let $\X$ be its minimal regular model. In the first step we prove that for every rigid vector bundle $E^{rig}\in\BB_{X^{rig}}^{\widehat{\X}}$ the algebraisation has a model $\E$ on $\X$ that satisfies $\E_{\op}\in\BB_{\X}$. This is slightly stronger than the statement that the algebraisation functor maps elements of the category $\BB_{X^{rig}}^{\widehat{\X}}$ to elements $E$ that satisfy $E_{\Cp}\in\BB_{X}$, because we specify the relevant model.
We remind of some notations that have been used before.\\
If $\Omega/\Gamma=X^{rig}$ is a Mumford curve, then the associated analytic reduction $\widetilde{X^{rig}}$ is a quotient $\widetilde{\Omega}/\Gamma$, and we have the commutative diagram
\[\xymatrix{\Omega \ar[r]^{u}\ar[d]^{R} & X^{rig} \ar[d]^{r} 
\\{\widetilde{\Omega} \ar[r]^{\widetilde{u}}} & {\widetilde{X^{rig}}} .}\]
For an algebraic vector bundle $E$ in the category $\BB_{X}^{\X}$ the vector bundle $[V^{\circ}_{\rho}\otimes r_{*}u_{*}\Ooo]^{\Gamma}$ is the formal completion of an algebraic model $\E$ of $E$ if $\Gamma$ acts by the representation associated to $E^{rig}$; this is proved in proposition \ref{modelandrep}. The lemma \ref{formalpullback} gives a criterion, whether the reduction modulo $p^{n}$ of a finite lift of $\E$ to another Mumford curve is trivial.
As a preparation we prove the following lemma.\\
Let $N\subset\Gamma$ be a cofinite normal subgroup, \[\xymatrix{{\Omega\ar[r]^{v}} & {\Omega/N}} \text{and }
\xymatrix{{\Omega\ar[r]^{u}} & {\Omega/\Gamma}}\] two Mumford curves and 
let \[ \xymatrix{{\Omega/N \ar[r]^{w}}& {\Omega/\Gamma}}\] be the canonical morphism such that $u=w\circ v$. 
\begin{lemma}
Let $E^{rig}:=[V\otimes u_{*}\Loc_{\Omega}]^{\Gamma}$ be a semistable vector bundle of degree 0 on $\Omega/\Gamma$ associated to $\rho:\Gamma\ra\Gl(V)$. Then $w^{*}E^{rig}$ is a semistable bundle of degree 0 on $\Omega/N$, and its associated representation of the Schottky group $\rho':N\ra\Gl(V')$ is isomorphic to $\rho_{|N}$.
\end{lemma}
\begin{proof}
Because $w$ is a finite morphism, the vector bundle $w^{*}E^{rig}$ is semistable of degree 0 and indeed is isomorphic to 
\[E'^{rig}:=[V'\otimes u_{*}\Loc_{\Omega}]^{N}\] 
for some representation $\rho':N\ra\Gl(V')$. The $\K$-vector space $V$ is a subspace of $u^{*}E^{rig}(\Omega)$ and the action of $\Gamma$ on $V$ is induced from the natural action of $\Gamma$ on $u^{*}E^{rig}(\Omega)$. Similarly the action on $V'\subset v^{*}E'^{rig}(\Omega)$ is induced by the natural action of $N$ on $v^{*}E'^{rig}(\Omega)$. The vector bundles $v^{*}E'^{rig}=v^{*}(w^{*}E^{rig})$ and $u^{*}E^{rig}$ are naturally isomorphic. Hence for all $\gamma\in N$ the diagram
\[\xymatrix{v^{*}E'^{rig} \ar[r]\ar[d] & \gamma^{*}v^{*}E'^{rig}\ar[d]\\
u^{*}E^{rig} \ar[r] & \gamma^{*}u^{*}E^{rig}},\] 
in which all arrows are isomorphisms, is commutative. The actions of $N$ on $V$ and on $V'$ coincide. It follows that $\rho_{|N}$ is isomorphic to $\rho'$.
\comment{ 
Because $u, v, w$ are topological coverings they are local isomorphism and there is an affinoid subdomain $V=\Sp(v_{*}\Loc_{\Omega}^{N}(V))$ of $\Omega/N$ and an affinoid subdomain $U=\Sp(u_{*}\Loc_{\Omega}^{\Gamma}(U))$ of $\Omega/\Gamma$ such that $w_{|V}:V\tilde{\ra} U$. 
Hence $v_{*}\Loc_{\Omega}^{N}(V)$ and $u_{*}\Loc_{\Omega}^{\Gamma}(U)$ are isomorphic affinoid algebras.\\
Since $E(\rho)$ is coherent it follows 
\[E(\rho)_{|U}=M\otimes u_{*}\Loc_{\Omega |U}^{\Gamma}=:\widetilde{M},\] $M$ is an $u_{*}\Loc_{\Omega}^{\Gamma}(U)$-module (by Kiehls theorem), also by Kiehls theorem it is $E(\rho')=\widetilde{M'}$ for a $v_{*}\Loc_{\Omega}^{N}(V)$-module $M'$. It follows \[w^{*}E(\rho)_{|U}=\widetilde{M\otimes_{u^{*}\Loc_{\Omega}^{\Gamma}(U)}v_{*}\Loc_{\Omega}^{N}(V)}\cong\widetilde{M}\] (since $u_{*}\Loc_{\Omega}^{\Gamma}(U)\cong v_{*}\Loc_{\Omega}^{N}(V)$).
Because of 
$w^{*}E(\rho)_{|U}\cong E(\rho')_{|U}$ it is $M\cong M'$\\
It is \[E(\rho)(U)=M\otimes u_{*}\Omega^{\Gamma}(U)=[V_{\rho}\otimes u_{*}\Loc_{\Omega}(U)]^{\Gamma}\] and \[E(\rho')(V)=M'\otimes v_{*}\Loc_{\Omega}^{N}(V)=[V_{\rho'}\otimes v_{*}\Loc_{\Omega}(V)]^{N}\cong [V_{\rho'}\otimes u_{*}\Loc_{\Omega}(U)]^{N}.\] It follows that $[V_{\rho}\otimes u_{*}\Loc_{\Omega}(U)]^{\Gamma}$ and $[V_{\rho'}\otimes u_{*}\Loc_{\Omega}(U)]^{N}$ are isomorphic $u_{*}\Loc_{\Omega}(U)^{\Gamma}$-modules: 
\[\xymatrix{{\phi:[V_{\rho}\otimes u_{*}\Loc_{\Omega}(U)]^{\Gamma} \ar[r]^{\cong}} & {[V_{\rho'}\otimes u_{*}\Loc_{\Omega}(U)]^{N}}}.\]
Because these modules are free an isomorphism between them is a square-matrix with elements in $u_{*}\Loc_{\Omega}(U)^{\Gamma}$.\\
Given an element $\vec{f}\in [V_{\rho}\otimes u_{*}\Loc_{\Omega}(U)]^{\Gamma}$ then for every $\gamma\in N$ holds \[\phi(\gamma\vec{f})=\phi(\rho(\gamma)\vec{f}\circ\gamma^{-1})=\phi\vec{f}.\] Furthermore \[\gamma\phi\vec{f}=\rho'(\gamma)\phi\vec{f}\circ\gamma^{-1}=\phi\vec{f}.\] It follows $\phi\rho_{|N}=\rho'\phi$.
}
\end{proof}

\begin{lemma}\label{formalpullback}
Let $E^{rig}$ be the semistable vector bundle of degree 0 on $\Omega/\Gamma$ associated to the representation $\rho:\Gamma\ra\Gl_{\rk}(\Ok)$.
Let $\E$ be the formal model $[V^{\circ}_{\rho}\otimes_{\K^{\circ}} \widetilde{u}_{*}R_{*}\Ooo]^{\Gamma}$ of $E^{rig}$. If $\rho_{|N}\equiv 1 \mod p^{n}$, then it follows that $w^{*}_{n}\E_{n}\cong \left([(\widetilde{v}_{*}R_{*}\Ooo)^{\rk}]^{N}_{n}\right)$. This is the trivial bundle of rank $\rk$ on $(\Omega/N)^{form}_{n}$.
\end{lemma}

\begin{proof}
Let $w^{form}$ be the morphism of formal schemes associated to $w$. By the preceding lemma \[(w^{form})^{*}\E=[V_{\rho_{|N}}\otimes \widetilde{v}_{*}R_{*}\Ooo]^{N}.\] Because of $\rho_{|N}\equiv 1 \mod p^{n}$ it is $(w^{form})^{*}_{n}\E_{n}\cong \left([(\widetilde{v}_{*}R_{*}\Ooo)]^{N}_{n}\right)^{\rk}$.  
\end{proof}
Using this criterion we can prove that every vector bundle $E\in\BB_{X}^{\X}$ has a model $\E$ on $\X$ such that $\E_{\op}$ is in $\BB_{\X_{\op}}$. The trivialising coverings are constructed explicitly.

Let $X$ be a $\K$-Mumford curve and $\X$ its minimal regular $\Ok$-model. Let $\rho: \Gamma\ra\Gl(\Ok)$ be a representation and $n$ a natural number. Let $E^{rig}$ be the associated semistable vector bundle of degree 0 and rank $\rk$ on $X^{rig}$ with formal model 
\[\E^{form}:=[V^{\circ}_{\rho}\otimes r_{*}u_{*}\Ooo]^{\Gamma}.\]
\begin{proposition}\label{covering}
There is a semistable $\Ok$-curve $\Y$ and a proper,
finitely presented $\Ok$-morphism $v:\Y\ra\X$,
whose generic fibre $v_{\K}:Y\ra X$ is a finite and \'etale Galois covering, such that $v_{K}^{an}:Y^{an}\ra X^{an}$ is in $\Cov^{ftop}(X^{an})$ and \[v_{n}^{*}\E_{n}=\Loc^{\rk}_{\Y_{n}}.\]
\end{proposition}
\begin{proof}
The rigidification of the Mumford curve $X$ is $X^{rig}=\Omega/\Gamma$. The formal model $\E$ of $E^{rig}$ exists because of $\Bild\rho\subset\Gl_{\rk}(\Ok)$ (cf. proposition \ref{modelandrep}). By corollary \ref{formelvectgaga}, the bundle $\E^{form}$ is the formal completion of an algebraic vector bundle $\E$ on $\X$. By lemma \ref{modelgaga}, $\E$ is an algebraic $\X$-model of the algebraisation $E$ of $E^{rig}$.\\
It follows
\[\E_{n}(U) =\E^{form}(U)_{n}=[V^{\circ}_{\rho}\otimes r_{*}u_{*}\Ooo]^{\Gamma}_{n}.\]
Let $\rho_{n}:\Gamma\ra\Gl_{\rk}(\Ok/p^{n})$ be the reduction of $\rho$ modulo $p^{n}$. The group $\Gl_{\rk}(\Ok/p^{n})$ is finite, hence $\rho_{n}$ factorises through a finite quotient \[\rho_{n}:\Gamma/N\ra\Gl_{\rk}(\Ok/p^{n}).\] Here $N$ is a normal cofinite subgroup of $\Gamma$; define $G:=\Gamma/N$. By lemma \ref{mumfordcover} the group $N$ is a Schottky group, and \[Y^{rig}:=\Omega/N\] is a Mumford curve of genus $(g-1)\cdot |\Gamma/N|+1$.
Let $\Y$ be the minimal regular model of the algebraisation of $Y^{rig}$. The model $\Y$ is a semistable curve. The schematic quotient \[\X':=\Y/G\] exists because $\Y$ is projective. By formal completion along the special fibre one obtains by lemma \ref{formalquot} the identity \[\widehat{\X'}=\widehat{\Y}/G=(\Omega^{form}/N)/G=\Omega^{form}/\Gamma\] as formal schemes. Therefore $\widehat{\X'}$ is the formal scheme associated to the Mumford curve $X$. Because of the uniqueness of the $\Ok$-model attached to a specific reduction (cf. theorem \ref{curvesdvr}) we obtain $\X'=\X$.\\
We obtained a semistable projective and flat $\K^{\circ}$-curve $\Y$ with a morphism \[w:\Y\ra\X.\] Because $\X$ is also projective and the base scheme is Noetherian, the morphism $w:\Y\ra\X$ is proper and finitely presented. 
Because $\Ok\ra\K$ is a flat morphism of rings, the generic fibre $w_{\K}$ is again a quotient morphism \[w_{\K}:Y\ra Y/G.\] The morphism $w_{\K}^{an}$ is a finite topological covering, hence $Y\ra Y/G$ is finite \'etale by \cite[Proposition 3.3.11]{Berk1993}. Alternatively: As $\Gamma$ acts on $\Omega$ without fixpoints, $G$ acts on $Y$ without fixpoints. So $Y\ra Y/G$ has a trivial decomposition group, so it is \'etale as proved in \cite[V Corollaire 2.4]{SGA1}.\\
After formal completion of $w$ along the closed fibre we obtain a morphism of formal schemes 
\[w^{form}:(\Omega/N,r_{N},\widetilde{\Omega/N})\ra (\Omega/\Gamma,r,\widetilde{\Omega/\Gamma}),\] thus the morphism
{\small\[
\xymatrix{
u^{form}:(\Omega,R,\widetilde{\Omega})\ar[r]^{v^{form}}& (\Omega/N,r_{N},\widetilde{\Omega/N}) \ar[r]^{w^{form}}& (\Omega/\Gamma,r,\widetilde{\Omega/\Gamma}).
}\]}
In this situation we can use lemma \ref{formalpullback}  -- where we have the same notations -- and conclude $w^{*}_{n}\E_{n}=\Loc_{\Y_{n}}^{\rk}$.
\end{proof}

The modulo $p^{n}$-trivialising covering $\Y\ra\X$, which has been constructed in proposition \ref{covering}, is an object in $\Cov_{\X}^{fpp}$, but it is indeed an object of a smaller category, which was used in a previous work by Deninger--Werner \cite[section 5]{Deni2003} and which will be defined now.
\begin{defin}
Let $R$ be a valuation ring with quotient field $Q$ of characteristic 0 and $\X$ an $R$-scheme. Let $\T_{\X}$ be the category whose objects are finitely presented and proper $G$-equivariant $R$-morphisms $\pi:\Y\ra\X$, where $G$ is a finite abstract group that operates by left action $R$-linearly on $\Y$ and trivially on $\X$, such that the generic fibre $\pi_{\Q}$ is an \'etale $G$-torsor.\\
Let $G'$ be a finite abstract group and the $G'$-invariant morphism $\pi':\Y'\ra\X$ an element of $\T_{\X}$. A morphism from $\pi:\Y\ra\X$ to $\pi':\Y'\ra\X$ is a commutative triangle
\[
\xymatrix{\Y \ar[rr]^{\phi} \ar[rd] && \Y' \ar[ld]\\
&\X}\]
together with a group morphism $\gamma:G\ra G'$ such that $\phi\circ g=\gamma(g)\circ\phi$, for $g\in G$.
\end{defin}
\begin{defin}
There is an obvious forgetful functor $\T_{\X}\ra \Cov_{\X}$. The objects of the full subcategory $\T_{\X}^{good}$ are the objects of $\T_{\X}$ that are mapped to $\Cov_{\X}^{good}$ by the forgetful functor.
\end{defin}
\begin{nota}
The covering $\pi:\Y\ra\Y/G$ (which was constructed in proposition \ref{covering}) is an object of $\T_{\X}^{good}$.
\end{nota}

\begin{nota}
Up to isomorphism the DW-representation $\rho_{E}$ of a vector bundle $E\in\BB_{X_{\Cp}}$ does not depend on the chosen model $\E$ of $E$. Therefore it is enough to calculate the representation for the canonical model $\E$ of $E$, which was defined in \ref{canmodel}.
\end{nota}
\begin{lemma}\label{modnequal}
Let $n$ be a natural number, $X$ a Mumford curve over $\K$ with minimal regular model $\X/\Ok$ and let $\overline{x}\in X(\overline{\K})$ be a geometric point of $X$. Let $E$ be a vector bundle in $\BB_{X}^{\X}$ and denote its canonical model by $\E$. Then $E^{rig}$ is in $\BB_{X^{rig}}^{\widehat{\X}}$. The reduction modulo $p^{n}$ of the DW-representation attached to $\E$ factorises through $\pi^{ftop}(X^{an},\overline{x})$ and is isomorphic to the reduction modulo $p^{n}$ of the PR-representation attached to $E^{rig}$.
\end{lemma}
\begin{proof}
The fact that $E^{rig}$ is in $\BB_{X^{rig}}^{\widehat{\X}}$ is proved in proposition \ref{modelandrep}. We constructed a covering $\pi:\Y\ra\X$ in $\T^{good}_{\X}$ forcing $\pi_{n}^{*}\E_{n}$ to be trivial. Since we have a $\T^{good}_{\X}$-covering here, we construct the DW-representation as it is done on page 37 of \cite{Deni2004}. 
\comment{

In the article the construction is done after a base change by $\Cp$ resp. $\op$. The consequence is that the base scheme is not Noetherian any more. We can avoid this as follows:
\begin{enumerate}
\item Consider the base change functors $F_{1}$ and $F_{2}$ in the category of $\Ok$-schemes, where $F_{1}(\X):=(\X\otimes_{\Ok}\op)\otimes_{\op}\on$ and  $F_{2}(\X):=(\X\otimes_{\Ok} \Ok/p^{n})\otimes_{\Ok/p^{n}}\overline{\Zp}/p^{n}$. $F_{1}$ and $F_{2}$ are isomorphic functors.
\item An $\Ok$-point $y$ of $\Y$ over an $\Ok$-point of $\X$ is a commutative triangle:
\[\xymatrix{\Spec \Ok\ar[rd]_{x}\ar[r]^{y}& {\Y}\ar[d]\\ &\X}.\] Because $F_{1}$ and $F_{2}$ are isomorphic functors they map the diagram to diagrams that are canonical isomorphic (this means the objects are canonical isomorphic and the morphisms are compatible with this canonical isomorphism). So $F_{1}(x)$ and $F_{2}(x)$ are canonical isomorphic points, which we identify writing $\overline{x}_{n}$
\end{enumerate}

}
Because of properness, the geometric point $\overline{x}\in X(\Cp)$ defines an unique point $x_{\op}$ in $\X(\op)$. Choose a point $\overline{y}\in Y(\Cp)$ above $\overline{x}$. It defines a point $y_{\op}\in\Y(\op)$. It determines a morphism, which factorises over $\pi_{1}^{ftop}(X_{\Cp},\overline{x})$ by remark \ref{fundamentalfactor} as in the following diagram:
\[\xymatrix{{\pi^{alg}(X_{\Cp},\overline{x})\ar[r]^{\phi_{y_{i}}}}\ar@{->>}[d]& {\Gal_{\overline{X}}Y_{i}\ar[r]^{\phi_{\overline{y}}}\ar[d]^{\phi_{y'_{i}}}}&{\Gal_{X_{\Cp}}Y} \ar[r]^{=}& {\Aut_{X}^{op}Y \ar[r]^{=}} & G^{op}\\
{\pi^{ftop}(X^{an}_{\Cp},\overline{x})\ar[r]^{\phi_{y'_{i}}}}&{\Gal_{X_{\Cp}}Y'_{i}\ar[ru]_{\phi_{\overline{y'}}}}}\] 
\comment{
\[\xymatrix{\pi_{1}^{alg}(\overline{X},\overline{x}) \ar[rr]^{\phi_{y}} \ar[rd] && {\Aut_{X}^{op}Y=G^{op}}\\
& {\pi_{1}^{ftop}(X,x)}\ar[ru]^{\phi'_{y}}},\]
}
By remark \ref{fundextension} we have the equation 
\[\pi_{1}^{ftop}(X^{an}_{\Cp},\overline{x})=\pi_{1}^{ftop}(X^{an},\overline{x}).\]
The reduction $\rho_{n}^{DW}$ is the morphism
\[\xymatrix{{\rho^{DW}_{n}:G^{op}\ar[r]}& {\Aut \E_{x_{n}}}\text{ ,  }\sigma\mapsto (y_{n}^{*})^{-1}\sigma^{*}y_{n}^{*}}\]
as in the following diagram where we abbreviated  
\[\Gamma(\Y_{n},(w^{form})_{n}^{*}[(V_{\rho}^{\circ}\otimes r_{*}u_{*}\Ooo)]^{\Gamma}_{n})\]
by the term $H^{0}$:
\[\xymatrix{{\sigma \ar @{|->}[r]} & ({\E_{x_{n}}\ar[r]^{(y_{n}^{*})^{-1}}}& H^{0}\ar[r]^{\sigma^{*}}&H^{0}\ar[r]^{y_{n}^{*}}& {\E_{x_{n}}})}
\]
The morphism $\sigma^{*}$ is defined as follows:
\[
\xymatrix{{H^{0}\ni f \ar @{|->}[r]}& {f\circ\sigma=\rho_{n}^{PR}(\sigma)f\in H^{0}}}.
\]
Therefore \[\rho^{DW}_{n}(\sigma)=(y_{n}^{*})^{-1}(\rho_{n}^{PR})(\sigma)(y_{n}^{*}).\] Hence $\rho^{DW}_{n}$ and $\rho_{n}^{PR}$ are isomorphic representations. Note that $\rho_{n}^{PR}$ is an element of $\Gl(\K^{\circ}/p^{n})$.
\end{proof}

\subsubsection{The general case}\text{ }\\
In this paragraph we compare the two representations for vector bundles that are only defined over $\Cp$. 
\begin{lemma}\label{generalcase}
Let $\K$ be a local number field and $n$ a natural number. Let $X$ be a Mumford curve over $\K$ with minimal regular model $\X/\Ok$. If $E$ is a vector bundle in $\BB_{X_{\Cp}}^{\X_{\op}}$ with canonical model $\E$, then $E^{rig}$ is in $\BB_{X_{\Cp}^{rig}}^{\widehat{\X_{\op}}}$. The reduction modulo $p^{n}$ of the DW-representation attached to $\E$ factorises through $\pi^{ftop}(X^{an},\overline{x})$ and is isomorphic to the reduction modulo $p^{n}$ of the PR-representation attached to $E^{rig}$.
\end{lemma}
\begin{proof}
The first assertion has been proved in proposition \ref{modelandrep}. By lemma \ref{reduction} there is a finite extension $\L^{(n)}$ of $\K$ and a vector bundle $\E^{(n)}$ on $\X_{\Loc_{\L^{(n)}}}$ such that \[\E^{(n)}\otimes \Loc_{\L^{(n)}}/p^{n}=\E\otimes\op/p^{n}.\] 
Now we can apply lemma \ref{modnequal} to see that there is an automorphism $A_{n}$ in $\Gl_{\rk}(\Loc_{\L^{(n)}})$ such that \[A_{n}^{-1}\rho^{PR}_{n}(E^{rig})A_{n}=\rho^{DW}_{n}(E).\]
\end{proof}
\begin{defin}
With the above notations define \[\rho^{cPR}(E^{rig}):=\lim_{n}\rho^{PR}_{n}(E^{rig}):\pi^{ftop}_{1}(X^{an},\overline{x})\ra\Gl_{\rk}(\K^{\circ})\] and call it {\em profinitely completed van der Put--Reversat representation}. By construction, it is continuous.
\end{defin}
\begin{theorem}\label{mainth}
Let $\K$ be local number field and let $\L$ be a complete subfield of $\Cp$ that is an extension of $\K$. Let $X$ be a Mumford curve over $\K$ with minimal regular model $\X/\Ok$ and $E$ a vector bundle in $\BB_{X_{\L}}^{\X_{\L^{\circ}}}$. Then the profinitely completed PR-representation extended to $\Cp$-vector spaces is isomorphic to the $DW$-representation. 
\end{theorem}
\begin{proof}
Let $\E$ be the canonical model of $E$. By the lemmas \ref{modnequal} and \ref{generalcase} there are projective systems $(A_{n})_{n\geq1}$, $(\rho^{PR}_{n})_{n\geq1}$ and $(\rho^{DW}_{n})_{n\geq1}$ such that \[A_{n}^{-1}\rho^{PR}_{n}(E^{rig})A_{n}=\rho^{DW}_{n}(E).\]
By functoriality of $\plim$ it is 
\[\plim\rho^{DW}_{n}(E^{rig})=\plim (A_{n}^{-1}\rho^{PR}_{n}(E)A_{n})=\plim (A_{n})^{-1} \plim\rho^{PR}_{n}(E) \plim A_{n}.\]
Furthermore with the notations before the theorem \ref{dewecp} it is \[\rho^{DW}(E)=\psi^{-1}_{\overline{x}}(\rho_{\E}(\gamma)\otimes\Cp)\psi_{\overline{x}}.\]
Therefore the representation $\rho^{DW}(E)$ is isomorphic to the profinitely completed PR-representation $\rho^{cPR}(E)\otimes\Cp$.
\end{proof}
\begin{corollar}
With the notations above, the representation $\rho^{DW}_{E}$ is isomorphic to a representation in the vector space $\L^{\rk}$.\\
Restricted to bundles in $\BB_{X^{rig}}^{\widehat{\X}}$ the profinitely completed PR-representation commutes with duals and tensor products, and the (original) RP-representation commutes with duals and tensor products as well.
\end{corollar}
The following lemma shows that representations of $\pi^{top}_{1}(X^{an},\overline{x})$ and continuous representations of $\pi^{ftop}_{1}(X^{an},\overline{x})$ are the same.
\begin{lemma}
\[\Hom(\Z^{g},\Gl_{\rk}(\K^{\circ}))\cong\Hom_{cont}(\widehat{\Z}^{g},\Gl_{\rk}(\K^{\circ}))\]
\end{lemma}
\begin{proof}
Every morphism $\phi$ in $\Hom(\Z^{g},\Gl_{\rk}(\K^{\circ}))$ is continuous if $\Z$ is endowed with the pro-finite topology and $\K^{\circ}$ carries the usual $p$-adic topology. This follows because $\K^{\circ}/p^{n}$ is torsion and $\Z^{g}$ is finitely generated. The group $\widehat{\Z}^{g}$ is the pro-finite completion of $\Z^{g}$, and the morphism $\phi$ has a unique completion $\widehat{\phi}\in\Hom_{cont}(\widehat{\Z}^{g},\Gl_{\rk}(\K^{\circ}))$.\\
Vice versa every $\psi\in\Hom_{cont}(\widehat{\Z}^{g},\Gl_{\rk}(\K^{\circ}))$ induces a morphism \[\Z^{g}\inj\widehat{\Z}^{g}\ra\Gl_{\rk}(\K^{\circ}).\]
\end{proof}

\subsubsection{Conclusion}\text{ }\\
In this thesis we proved the following: 
\begin{enumerate}
\item The categories $\BB_{X}^{\X}$ and $\BB_{X^{rig}}^{\widehat{\X}}$ are equivalent. The equivalence is induced by rigidification of vector bundles (cf. proposition \ref{modelandrep}). 
\item The DW-representation attached to a vector bundle in $\BB_{X}^{\X}$ factorises over the finite topological fundamental group $\pi^{ftop}_{1}(X^{an},\overline{x})$ of $X$ and is isomorphic to the profinitely completed PR-representation attached to the rigidification of $E$ (cf. theorem \ref{mainth}).
\item The category $\BB_{X}^{\X}$ is equivalent to the category of continuous $\K$-vector space-representations of $\pi^{ftop}(X^{an},\overline{x})$ resp. to the category of $\K$-vector space-representations of $\pi^{top}(X^{an},\overline{x})$ for which a basis can be found such that the associated matrices are in $\Gl_{\rk}(\K^{\circ})$. The equivalence commutes with duals, tensor products, extensions and internal homs.
\end{enumerate}
\section{Illustrations}
In this last section we show various illustrations. In the first subsection we have a look on the different fundamental groups of a Tate curve. In the second and third subsection we make the considered categories of representations more explicit in the cases of Tate curves and Mumford curves of genus 2.
\subsection{Fundamental groups of a Tate curve}
Let $X$ be a Tate curve over a local number field $\K$. Its rigidification $X^{rig}$ is a quotient $\GmK/\Gamma$ where $\Gamma$ is the cyclic group that is generated by $\begin{pmatrix} q & 0 \\ 0 & 1 \end{pmatrix}$ for an element $q\in\K$ that satisfies $0<\n{q}< 1$.
The algebraic fundamental group of $X$ is \[\pi_{1}^{alg}(X_{\Cp},1)=(\widehat{\Z}^{2})^{ab}.\] Because elliptic curves are abelian varieties, a cofinite system of its Galois coverings is the system of $N$-multiplications. Such a covering is not topological, but it factorises as
\[\xymatrix{\GmK/\Gamma\ar[d]^{z\mapsto z^{n}}\\
\GmK/\langle q^{n}\rangle\ar[d]\\
\GmK/\Gamma}.\]
The second morphism is a topological covering, but not the first one. We remind the reader that the different fundamental groups were defined in definitions \ref{covtypes} and \ref{funddef}.
Yves Andr\'e shows in paragraph \cite[II 2.3.2]{Andr2003a} that 
\[\pi_{1}^{temp}(X_{\Cp}^{an},1)=\pi_{1}^{alg}(\mathbb{G}_{m,\Cp},1)\times\pi_{1}^{top}(X_{\Cp}^{an},1).\]
The profinite completion of $\pi_{1}^{temp}(X_{\Cp}^{an},1)$ is $\pi_{1}^{alg}(X_{\Cp}^{an},1)$, therefore we obtain by profinite completion:
\[\pi_{1}^{alg}(X_{\Cp},1)=\pi_{1}^{alg}(\mathbb{G}_{m,\Cp},1)\times\pi_{1}^{ftop}(X_{\Cp}^{an},1)=\widehat{\Z}(1)\times\widehat{\Z}.\]
This result was already indicated by the above factorisation of the $N$-multiplication. 
\subsection{Vector bundles on a Tate curve}  
Over an algebraically closed field a representation $\Z\ra\Gl_{\rk}(\Cp)$ is, up to isomorphism, uniquely determined by its Jordan normal form, up to permutation of Jordan blocks. It follows that objects in the skeleton of $\BB_{X_{\Cp}}^{\X_{\op}}$ correspond to $\phi$-bounded matrices over $\Cp$ in Jordan normal form up to block permutation.\\
\begin{ex}
The representation associated to a vector bundle of rank 2 in $\BB_{X_{\Cp}}^{\X_{\op}}$ is isomorphic to a representation that is represented by one of the following matrices:
\[ \left( \begin{array}{cc}
x & 0  \\
0 & y \end{array} \right),
\left( \begin{array}{cc}
x & 0  \\
0 & x \end{array} \right),
\left( \begin{array}{cc}
x & 1  \\
0 & x \end{array} \right)  
 \]
  with the conditions $x,y\in\op^{\times}$ and $x\ne y$. Here remind that two matrices describe isomorphic representations if and only if they are conjugated via a matrix in $\Gl_{2}(\Cp)$; e.g. 
\[ \left( \begin{array}{cc}
1 & 0  \\
0 & p \end{array} \right)  
\left(\begin{array}{cc}
x & p  \\
0 & x \end{array} \right)
\left( \begin{array}{cc}
1 & 0  \\
0 & p^{-1} \end{array} \right)
=
\left(\begin{array}{cc}
x & 1  \\
0 & x \end{array} \right).\]
\end{ex}
Because the Jordan normal form is only useful over an algebraically closed field, we remind here of the so-called {\em rational normal form}. Let $Q_{i}(X)=\sum_{j=0}^{s} a_{j}X^{j}$ be a polynomial in $\K[X]$, we associate to $Q(X)$ a {\em companion matrix} 
\[C_{Q_{i}}:=\left( \begin{array}{cccccc}
0 & 0 & 0 &\cdots & 0 & -a_{0}  \\
1 & 0 & 0 &\cdots & 0 & -a_{1}  \\
0 & 1 & 0 &\cdots & 0 & -a_{2}  \\
\vdots & \vdots & \vdots &\ddots & \vdots & \vdots  \\
0 & 0 & 0 &\cdots & 1 & -a_{s-1}
\end{array} \right).\]
With this definition it is true that every matrix $M\in\M_{n,n}(\K)$ is conjugate over $\K$ to a matrix in rational canonical form; this is a matrix of the form
\[\left( \begin{array}{cccc}
C_{Q_{1}} & 0 &\cdots & 0 \\
0 & C_{Q_{2}} & \cdots & 0  \\
\vdots & \vdots &\ddots & \vdots  \\
0 & 0 &\cdots & C_{Q_{r}}
\end{array} \right),\]
where $\prod_{i=1}^{r}Q_{i}$ is a factorisation of the characteristic polynomial of $M$ in monic polynomials, such that $Q_{r}$ is the minimal polynomial of $M$ and $Q_{i}$ divides $Q_{i+1}$ for $i=1,\ldots,r-1$. The rational canonical form is unique.
\begin{lemma}
Let $\K$ be a subfield of $\Cp$. A matrix $M\in\M_{n,n}(\K)$ is conjugate over $\K$ to a matrix $M'\in\M_{n,n}(\Ok)$ if and only if the characteristic polynomial of $M$ is in $\Ok[X]$. 
\end{lemma}
\begin{proof}
Assume $M$ is conjugate to $M'\in\M_{n,n}(\Ok)$. Because the characteristic polynomial is invariant under conjugation, it can be calculated from $M'$. In the calculation only multiplication, addition and subtraction is involved, hence the polynomial is in $\Ok[X]$.\\
Assume that the characteristic polynomial $P(X)$ of $M$ is in $\K^{\circ}[X]$. The ring $\K[X]$ is factorial. The prime factors of $P(X)$ are polynomials in $\K^{\circ}[X]$, since the algebraic closure of $\K$ has an integrally closed ring of integers. Therefore in the decomposition $\prod_{i}Q_{i}(X)$ of $P(X)$ all the $Q_{i}(X)$ are elements of $\Ok[X]$.
\comment{
We prove that every zero of $P(X)$ is in $\op$. Assume the contrary and let $\alpha\in\Cp$ be a zero of $P(X)$ that is not an element of $\op$. Let $p^{-q}=|\alpha|_{\Cp}$, by assumption $q>0$. Then 
\[\begin{array}{lll}
0&=&P(\alpha)=\sum_{i=0}^{n} a_{i}\alpha^{i}\\
0&=&\sum_{i=0}^{n} p^{qn}a_{i}\alpha^{i}\\
0&=&\sum_{i=0}^{n} p^{q(n-i)}a_{i}(\alpha\cdot p^{q})^{i}.\\
\end{array}\]
Let $\m$ be the maximal ideal in $\op$. The first summand of 
\[\sum_{i=0}^{n} p^{q(n-i)}a_{i}(\alpha\cdot p^{q})^{i}\] is a unit in $\op$ and all other summands are in $\m$. Because the sum is zero it is also in $\m$. This is a contradiction. It follows that all zeros are in $\op$.} 
Therefore the companion matrices in the rational canonical form have integral entries.
\end{proof}
\begin{corollar}
Every semistable vector bundle of degree zero on the Tate curve is uniquely characterised by the canonical rational form of its associated matrix. It is an object in $\BB_{X_{\K}}^{\X_{\K^{\circ}}}$ if and only if this rational form is a matrix in $\Gl(\K^{\circ})$. 
\end{corollar}

\subsection{Vector bundles on a Mumford curve of genus 2}
Let $\K$ be a local number field, and let $X$ be a $\K$-Mumford curve of genus 2. In this case an isomorphism class of semistable vector bundles of degree zero on $X_{\Cp}$ is associated to an isomorphism class of $\phi$-bounded representations \[\rho:\Z^{2}\ra\Gl_{\rk}(\Cp).\]
An isomorphism class of $\Cp$-vector space respresentations of the free group generated by $\gamma_{1}, \gamma_{2}$ is determined by the pair of matrices \[(\rho(\gamma_{1}),\rho(\gamma_{2}))\in\Gl_{\rk}(\Cp)^{2}\] up to simultanuous conjugation. It is a very difficult problem to classify all conjugacy classes of these pairs. 
As an example we discuss the easiest case, which is the one of rank 2 matrices over the algebraically closed field $\Cp$.\\
Let $\mathcal{M}'_{2,2}$ be the set of matrices in $\mathcal{M}_{2,2}(\Cp)$ that have trace zero. We have a direct sum decomposition $\mathcal{M}_{2,2}=\Cp\oplus\mathcal{M}'_{2,2}$.
By \cite[2.4 Example]{Kraf1996} there is a Zariski-dense subset $U'$ of $(\mathcal{M}'_{2,2})^{2}$ where every pair $(A',B')\in U'$ is conjugate to a pair of the form
\[ \left(\left( \begin{array}{cc}
t & 0  \\
0 & -t \end{array} \right),
\left( \begin{array}{cc}
s & 1  \\
c & -s \end{array} \right)\right) \text{, with } t,c\ne 0\text{ and}
\]
\[t^2=\frac{1}{2}\tr A'^2 \text{, } s^2=\frac{(\tr A'B')^{2}}{\tr A'^{2}} \text{ and } c=\frac{1}{2}\tr B'^{2}-s^2.\] This classifies the pairs $(A,B)$ of a Zariski-dense subset $U\in\mathcal{M}_{2,2}^{2}$ up to simultaneous conjugation by adding the respective traces. The set \[\{(A,B)\in U \suchthat (\tr A^2-\tr A)(\tr B^{2}-\tr B)\ne 0\}\] is Zariski-dense in $\Gl_{2}(\Cp)^{2}$. Hence we classified a Zariski-dense subset of $\Gl_{2}(\Cp)^{2}$ up to conjugation.  
\bibliography{lit}

\providecommand{\bysame}{\leavevmode\hbox to3em{\hrulefill}\thinspace}
\begin{thebibliography}{GvdP80}

\bibitem[AM69]{Atiy1969}
M.~F. Atiyah and I.~G. Macdonald, \emph{Introduction to commutative algebra},
  Addison-Wesley Publishing Co., Reading, Mass.-London-Don Mills, Ont., 1969.

\bibitem[And03a]{Andr2003b}
Yves Andr{\'e}, \emph{On a geometric description of {${\rm Gal}(\overline{\bf
  Q}\sb p/{\bf Q}\sb p)$} and a {$p$}-adic avatar of {$\widehat{GT}$}}, Duke
  Math. J. \textbf{119}:1 (2003), 1--39.

\bibitem[And03b]{Andr2003a}
Yves Andr{\'e}, \emph{Period mappings and differential equations. {F}rom
  {$\mathbb{C}$} to {$\mathbb{C}_{p}$}}, MSJ Memoirs, vol.~12, Mathematical
  Society of Japan, Tokyo, 2003, Tohoku-Hokkaido lectures in arithmetic
  geometry.

\bibitem[Ber90]{Berk1990}
Vladimir~G. Berkovich, \emph{Spectral theory and analytic geometry over
  non-{A}rchimedean fields}, Mathematical Surveys and Monographs, vol.~33,
  American Mathematical Society, Providence, RI, 1990.

\bibitem[Ber93]{Berk1993}
Vladimir~G. Berkovich, \emph{\'{E}tale cohomology for non-{A}rchimedean
  analytic spaces}, Inst. Hautes \'Etudes Sci. Publ. Math.:78 (1993), 5--161
  (1994).

\bibitem[Ber99]{Berk1999}
Vladimir~G. Berkovich, \emph{Smooth {$p$}-adic analytic spaces are locally
  contractible}, Invent. Math. \textbf{137}:1 (1999), 1--84.

\bibitem[BGR84]{Bosc1984}
S.~Bosch, U.~G{\"u}ntzer, and R.~Remmert, \emph{Non-{A}rchimedean analysis},
  Grundlehren der Mathematischen Wissenschaften, vol. 261, Springer-Verlag,
  Berlin, 1984.

\bibitem[dJ95]{Jong1995}
A.~J. de~Jong, \emph{\'{E}tale fundamental groups of non-{A}rchimedean analytic
  spaces}, Compositio Math. \textbf{97}:1-2 (1995), 89--118.

\bibitem[DWa]{Deni2004b}
Christopher Deninger and Annette Werner, \emph{Line bundles and $p$-adic
  characters}, \mbox{math.AG/0407511}.

\bibitem[DWb]{Deni2003}
Christopher Deninger and Annette Werner, \emph{Vector bundles and p-adic
  representations I}, \mbox{math.NT/0309273}.

\bibitem[DWc]{Deni2004}
Christopher Deninger and Annette Werner, \emph{Vector bundles on p-adic curves
  and parallel transport}, \mbox{math.AG/0403516}.

\bibitem[EGA]{EGA}
A.~Grothendieck, \emph{\'{E}l\'ements de g\'eom\'etrie alg\'ebrique. {I - IV}},
  Inst. Hautes \'Etudes Sci. Publ. Math.:4,8,11,17,20,24,28,32 (1960---1967).

\bibitem[EGAn]{EGAn}
A.~Grothendieck and Jean~A. Dieudonn\'e, \emph{{\'El\'ements de g\'eom\'etrie
  alg\'ebrique. I}}, {Die Grundlehren der mathematischen Wissenschaften. 166.
  Berlin-Heidelberg-New York: Springer-Verlag. IX, 466 p. }, 1971.

\bibitem[Fal83]{Falt1983}
Gerd Faltings, \emph{Semistable vector bundles on {M}umford curves}, Invent.
  Math. \textbf{74}:2 (1983), 199--212.

\bibitem[Fal03]{Falt2003}
Gerd Faltings, \emph{A $p$-adic simpson correspondence}, Notes (2003).

\bibitem[FvdP04]{Fres2004}
Jean Fresnel and Marius van~der Put, \emph{Rigid analytic geometry and its
  applications}, Progress in Mathematics, vol. 218, Birkh\"auser Boston Inc.,
  Boston, MA, 2004.

\bibitem[Gir71]{Gira1971}
Jean Giraud, \emph{Cohomologie non ab\'elienne}, Springer-Verlag, Berlin, 1971,
  Die Grundlehren der mathematischen Wissenschaften, Band 179.

\bibitem[Gro57]{Grot1957}
A.~Grothendieck, \emph{Sur la classification des fibr\'es holomorphes sur la
  sph\`ere de {R}iemann}, Amer. J. Math. \textbf{79} (1957), 121--138.

\bibitem[GvdP80]{Gerr1980}
Lothar Gerritzen and Marius van~der Put, \emph{Schottky groups and {M}umford
  curves}, Lecture Notes in Mathematics, vol. 817, Springer, Berlin, 1980.

\bibitem[Iha66]{Ihar1966}
Yasutaka Ihara, \emph{On discrete subgroups of the two by two projective linear
  group over {${\p}$}-adic fields}, J. Math. Soc. Japan \textbf{18} (1966),
  219--235.

\bibitem[K{\"o}p74]{Kpf1974}
Ursula K{\"o}pf, \emph{\"{U}ber eigentliche {F}amilien algebraischer
  {V}ariet\"aten \"uber affinoiden {R}\"aumen}, Schr. Math. Inst. Univ.
  M\"unster (2):Heft 7 (1974), iv+72.

\bibitem[KP96]{Kraf1996}
Hanspeter Kraft and Claudio Procesi, \emph{Classical invariant theory, a
  primer}, \mbox{http://www.math.unibas.ch/~kraft/Papers/KP-Primer.pdf}.

\bibitem[Liu02]{Liu2002}
Qing Liu, \emph{Algebraic geometry and arithmetic curves}, Oxford Graduate
  Texts in Mathematics, vol.~6, Oxford University Press, Oxford, 2002.

\bibitem[L{\"u}t93]{Ltke1993}
W.~L{\"u}tkebohmert, \emph{Riemann's existence problem for a {$p$}-adic field},
  Invent. Math. \textbf{111}:2 (1993), 309--330.

\bibitem[Mum63]{Mumf1962}
David Mumford, \emph{Projective invariants of projective structures and
  applications}, Proc. Internat. Congr. Mathematicians (Stockholm, 1962), Inst.
  Mittag-Leffler, Djursholm, 1963, pp.~526--530.

\bibitem[Mum72]{Mumf1972}
David Mumford, \emph{An analytic construction of degenerating curves over
  complete local rings}, Compositio Math. \textbf{24} (1972), 129--174.

\bibitem[NS65]{Nara1965}
M.~S. Narasimhan and C.~S. Seshadri, \emph{Stable and unitary vector bundles on
  a compact {R}iemann surface}, Ann. of Math. (2) \textbf{82} (1965), 540--567.

\bibitem[Ray74]{Rayn1974}
Michel Raynaud, \emph{G\'eom\'etrie analytique rigide d'apr\`es {T}ate,
  {K}iehl,{$\cdots $}}, Table Ronde d'Analyse non archim\'edienne (Paris,
  1972), Soc. Math. France, Paris, 1974, pp.~319--327. Bull. Soc. Math. France,
  M\'em. No. 39--40.

\bibitem[SGA1]{SGA1}
A.~Grothendieck, \emph{Rev\^etements \'etales et groupe fondamental ({SGA} 1)},
  Documents Math\'ematiques (Paris), 3, Soci\'et\'e Math\'ematique de France,
  Paris, 2003, S\'eminaire de g\'eom\'etrie alg\'ebrique du Bois Marie
  1960--61. Lecture Notes in Math., 224, Springer, Berlin.

\bibitem[vdPR86]{Put1986}
Marius van~der Put and Marc Reversat, \emph{Fibr\'es vectoriels semi-stables
  sur une courbe de {M}umford}, Math. Ann. \textbf{273}:4 (1986), 573--600.

\bibitem[Wei38]{Weil1938}
A.~Weil, \emph{{G\'en\'eralisation des fonctions ab\'eliennes.}}, J. Math.
  Pures Appl. \textbf{IX. S\'er. 17} (1938), 47--87.

\end{thebibliography}
\end{document}